\RequirePackage{fix-cm}
\documentclass[smallextended]{svjour3}       
\smartqed  
\usepackage{graphicx}
\usepackage{amsfonts}
\usepackage{booktabs}
\usepackage[fleqn]{amsmath}
\usepackage{epsfig}
\usepackage{amsmath}
\usepackage{tikz}
\usepackage[latin1]{inputenc}
\usepackage{verbatim}
\usepackage{mathptmx}
\usepackage{graphicx}

\usepackage{amssymb}
\usepackage{mathrsfs}
\usepackage[colorlinks,citecolor=blue,urlcolor=blue]{hyperref}
\usepackage{listings,matlab-prettifier} 

\newtheorem{algorithm}{Algorithm}[section]

\begin{document}

\title{An overlapping domain decomposition splitting algorithm for stochastic nonlinear Schr\"odinger equation\thanks{The research of Lihai Ji is partly supported by the National Natural Science Foundation of China (12171047, 11971458).}
}

\titlerunning{ODDS algorithm for stochastic NLS equation}        

\author{Lihai Ji
}


\institute{L. Ji \at
              Institute of Applied Physics and Computational Mathematics, Beijing 100094, China \\
              \email{jilihai@lsec.cc.ac.cn}           
}

\date{Received: date / Accepted: date}

\maketitle

\begin{abstract}
A novel overlapping domain decomposition splitting algorithm based on a Crank--Nisolson method is developed for the stochastic nonlinear Schr\"odinger equation driven by a multiplicative noise with non-periodic boundary conditions. The proposed algorithm can significantly reduce the computational cost while maintaining the similar conservation laws. Numerical experiments are dedicated to illustrating the capability of the algorithm for different spatial dimensions, as well as the various initial conditions. In particular, we compare the performance of the overlapping domain decomposition splitting algorithm with the stochastic multi-symplectic method in [S. Jiang, L. Wang and J. Hong, Commun. Comput. Phys., 2013] and the finite difference splitting scheme in [J. Cui, J. Hong, Z. Liu and W. Zhou, J. Differ. Equ., 2019]. We observe that our proposed algorithm has excellent computational efficiency and is highly competitive. It provides a useful tool for solving stochastic partial differential equations.
\end{abstract}

\section{Introduction}
\label{1.1}
The main purpose of this work is to propose an innovative overlapping domain decomposition splitting (ODDS for short) algorithm for the stochastic nonlinear Schr\"odinger (NLS) equation driven by a multiplicative noise
\begin{equation}
\begin{cases}\label{NLS}
{\rm i}du=\big[\Delta u+\lambda|u|^{2}u\big]dt+\varepsilon u\circ dW(t),\quad t\in(0,T],\\[1.5mm]
u(0,x)=u_0(x),\quad x\in D\subset\mathbb{R}^{d},\quad d\geq1,
\end{cases}
\end{equation}
where $\varepsilon>0$, $\lambda\in\mathbb{R}$ and $W$ is an $L^2(D;\mathbb{R})$-valued $Q$-Wiener process defined on a complete filtered probability space $(\Omega,\mathcal{F},\{\mathcal{F}_t\}_{t\in[0,T]},\mathbb{P})$. More precisely, $W(t)$ has the following Karhunen--Lo\`eve expansion
\begin{equation*}
W(t)=\sum_{k\in\mathbb{N}^d}Q^{\frac12}e_k\beta_{k}(t),\quad t\in[0,T]
\end{equation*}
with$\{e_{k}\}_{k\in \mathbb{N}^d}$ being an orthonormal basis of $L^2(D;\mathbb{R})$, $\{\beta_k\}_{k\in{\mathbb N^d}}$ being a sequence of real-valued mutually independent and identically distributed Brownian motions, and $Qe_{k}=\eta_k e_k$ for $\eta_k\geq 0$, $k\in{\mathbb N}^d$. For convenience, we always consider the equivalent It\^o form of \eqref{NLS}
\begin{equation}\label{Ito}
{\rm i}du=\big[\Delta u+\lambda|u|^{2}u-\frac{{\rm i}}{2}\varepsilon^{2}F_{Q}u\big]dt+\varepsilon udW(t)
\end{equation}
with the initial value $u(0)=u_0$ and $F_{Q}:=\sum_{k\in\mathbb{N}^d}(Q^{\frac12}e_{k})^{2}$.

In the last two decades, much progress has been made in theoretical analysis and numerical approximation for the stochastic NLS equation. To numerically inherit the charge conservation law and the geometric structure of \eqref{NLS}, \cite{RefCH2016,RefCHJ2017,RefCHLZJ2017,JWH2013} propose the stochastic symplectic and multi-symplectic algorithms. Particularly, the authors in \cite{RefCHJS2022} applies the large deviation principle to investigate the probabilistic superiority of the stochastic symplectic algorithms. To preserve the ergodicity of the numerical solution of \eqref{NLS}, \cite{RefCHLZ2017,HW2019,HWZ2017} study the ergodic numerical approximations. To reduce the computational cost of \eqref{NLS}, \cite{HWZ2019} designs a parareal algorithm and \cite{CHLZ2019,RefLiu,RefLiu2} propose the splitting algorithm, respectively. For more details about other kinds of numerical approximations of the stochastic NLS equation, we refer to \cite{BDD2005,RefBouard3,RefBouard2,RefBouard1,CHP2016} and references therein. These existing semi-discretizations and full discretizations for the stochastic NLS equation mentioned above are all investigated under the assumption of homogeneous or periodic boundary conditions. It is worth to point out that the soliton solution of the nonlinear dispersive wave propagation problems in a very large or unbounded domain for the stochastic NLS  equation is an interesting and important subject in applications (see, e.g., \cite{BCIRG1995}). This motives us to construct high efficient and numerical stable algorithms for the $d$-dimensional stochastic NLS equation \eqref{NLS} in a large spatial domain with non-zero or non-periodic boundary conditions.

To this end, we first apply the splitting technique in \cite{RefLiu} to split the equation \eqref{NLS} and get a deterministic linear PDE and a nonlinear stochastic PDE:
\begin{align}
{\rm i}du^{[1]}&=\Delta u^{[1]}dt,\label{sub_1}\\[1.2mm]
{\rm i}du^{[2]}&=\lambda|u^{[2]}|^2|u^{[2]}|dt+\varepsilon u^{[2]}\circ dW(t).\label{sub_2}
\end{align}
Then, for the subsystem \eqref{sub_2} we can get the analytic solution due to the point-wise conservation law $|u(t,x)|=|u_{0}(x)|$. The key issues lie in the numerical approximation for the subsystem \eqref{sub_1}, we first discretize it based on the Crank--Nicolson scheme in the temporal direction and get a temporal semi-discretization. 

To overcome the difficulties introduced by the non-periodic boundary conditions, we use the Chebyshev pseudo-spectral interpolation idea in space. To efficiently exploit modern high performance computing platforms, it is essential to design high performance algorithms. Domain decomposition method provides a useful tool to develop fast and efficient solvers for stochastic PDEs with a large number of random inputs. The non-overlapping domain decomposition method for PDEs with random coefficients is first proposed in \cite{SBG2009} and then extended by \cite{SS2014} to quantify uncertainty in large-scale simulations. We refer to \cite{DKPPS2018,LTT2010,TST2017,TST2018} and references therein for more details about the theory and applications of the domain decomposition method to PDEs with the random input.

We combine the Chebyshev interpolation idea and the overlapping domain decomposition method to approximate the temporal semi-discretization of \eqref{sub_1} and thus obtain a full discretization of \eqref{sub_1}. The explicitness of the solution of \eqref{sub_2} together with the full discretization of \eqref{sub_1} gives us an ODDS algorithm for the stochastic NLS equation \eqref{NLS}. Finally, several numerical examples for the stochastic NLS equation in one and two-dimensional spaces are presented to illustrate the capability of the proposed algorithm, which can be calculated high efficient. To the best of our knowledge, this is the first domain decomposition result of numerical approximations for stochastic PDEs whose the stochasticity comes from the stochastic source.

The rest of this paper is organized as follows. In Section \ref{ODDS}, we present and analyze the ODDS algorithm for the stochastic NLS equation. In Section \ref{2.1}, we show the algorithm for the one-dimensional stochastic NLS equation. In Section \ref{2.2}, we focus on studying the ODDS algorithm for the two-dimensional case. Section \ref{ne} contains some numerical experiments for the stochastic NLS equation to demonstrate the accuracy and efficiency of the proposed algorithm. Concluding remarks are given in Section \ref{cr}.
\section{The ODDS algorithm for the stochastic NLS equation}\label{ODDS}
In this section, we devote to obtaining the ODDS algorithm for the stochastic NLS equation in one and multi-dimensional spaces. 
\subsection{An ODDS algorithm for the one-dimensional stochastic NLS equation}\label{2.1}
This part concentrates mainly on demonstrating an ODDS algorithm for the following stochastic nonlinear problem:
\begin{equation}\label{NLS_one}
{\rm i}du=\big[u_{xx}+\lambda|u|^{2}u\big]dt+\varepsilon u\circ dW(t),\quad t\in(0,T]
\end{equation}
with an initial datum 
$$u(0,x)=u_0(x),\quad x\in D=[x_L,x_R]$$ 
and the boundary conditions 
\begin{equation}\label{bc}
u(t,x_L)=f(t),\quad u(t,x_R)=g(t),\quad t\in(0,T].
\end{equation}

As is well known, the stochastic NLS equation \eqref{NLS} possesses the charge conservation law under the homogeneous or periodic boundary conditions (see, e.g., \cite[Proposition 4.4]{RefBouard0}), that is
\begin{equation}\label{7}
\int_D|u(t,x)|^2dx=\int_D|u_0(x)|^2dx,\quad\mathbb{P}\text{-}a.s.
\end{equation}
for all $t\in[0,T]$. Furthermore, if we define the Hamiltonian
\begin{equation*}
H(u)=\frac{1}{2}\int_{D}|\nabla u(t,x)|^{2}{\rm d}x-\frac{1}{4}\int_{D}|u(t,x)|^{4}{\rm d}x,
\end{equation*}
then the averaged energy $\mathbb{E}[H(u(t))]$ satisfies (see, e.g., \cite[Proposition 4.5]{RefBouard0})
\begin{equation}\label{energy}
\begin{split}
\mathbb{E}[H(u(t))]=\mathbb{E}[H(u_{0})]+\frac{\varepsilon^{2}}{2}\int_{0}^{t}\int_{D}|u(s,x)|^{2}\sum_{k\in\mathbb{N}^d}\big|\nabla(Q^{\frac12} e_{k}(x))\big|^{2}{\rm d}x{\rm d}s.
\end{split}
\end{equation}
In general, there are no charge conservation law and the averaged energy evolution law for the stochastic NLS equation with the boundary conditions \eqref{bc}.

\vspace{3mm}
\noindent{\bf (a). Operator splitting}
\vspace{2mm}

We use a splitting technique which was proposed in \cite{RefLiu} to discretize \eqref{NLS_one}. Denote 
\begin{align}
{\rm i}du&=L(u)dt,\label{NLS_two1}\\[1.2mm]
{\rm i}du&=N(u)dt,\label{NLS_two2}
\end{align}
where 
\begin{align*}
L(u)=u_{xx},\quad N(u)=\lambda|u|^2u-\frac{{\rm i}}{2}\varepsilon^2F_Qu+\varepsilon u\dot{W}.
\end{align*}
Note that for the nonlinear subsystem \eqref{NLS_two2}, we have the following useful result. We refer readers to \cite[Proposition 3.1]{RefLiu} for more details.
\begin{proposition}\label{pro1}
	Assume the initial datum $u_{0}(x)$ is $\mathcal{F}_{0}$-measurable $L^{2}$-valued random variable. Then the solution of 
	\begin{equation*}
	{\rm i}du=N(u)dt
	\end{equation*}
	is given by
	\begin{equation*}
	u(t,x)=u_{0}(x)\exp\big\{-{\rm i}\big(t\lambda|u_{0}(x)|^{2}+\varepsilon W(t,x)\big)\big\}.
	\end{equation*}
	Specially, $|u(t,x)|=|u_{0}(x)|$.
\end{proposition}
Motivated by this proposition, we get the following recursion in $T_n=(t_n,t_{n+1}]$, $t_n=n\tau$, $n\in\{0,1,\ldots,N-1\}$:
\begin{equation}\label{11}
u^{n+1}=\exp\big\{-{\rm i}\big(\tau\lambda|u^{n}|^{2}+\varepsilon \Delta W^{n+1}\big)\big\}u^{n},
\end{equation}
with $u^0=u_0$ and $\Delta W^{n+1}:=W(t_{n+1},x)-W(t_n,x)$.

\vspace{3mm}
\noindent{\bf (b). Overlapping domain decomposition method}
\vspace{2mm}

Now we are in a position to approximate the deterministic linear subsystem \eqref{NLS_two1}. Denoting by $p$ and $q$ the real and imaginary parts of the solution $u$ of \eqref{NLS_two1}, which satisfy 
\begin{align}
dp=q_{xx}dt,\quad dq=-p_{xx}dt.
\end{align}
Applying the Crank--Nicolson method to discretize the above equations in the temporal direction yields
\begin{equation}\label{first_order_system}
\begin{split}
p^{n+1}=p^{n}+\frac{\tau}{2}\Big(q_{xx}^{n+1}+q_{xx}^{n}\Big),\quad q^{n+1}=q^{n}-\frac{\tau}{2}\Big(p_{xx}^{n+1}+p_{xx}^{n}\Big)
\end{split}
\end{equation}
for all $n=0,1,\ldots,N-1$

\begin{figure}[th!]
	\begin{center}
		\includegraphics[width=1\textwidth]{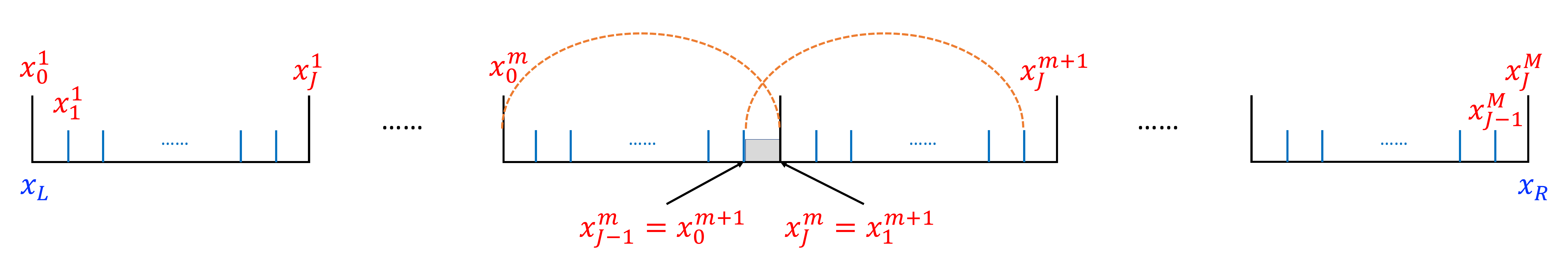}
		\caption{Basic idea of the overlapping domain decomposition method.}\label{oddM_fig}
	\end{center}
\end{figure}

Before we come to the spatial discretization of \eqref{first_order_system}, let us introduce some basic concepts of the overlapping domain decomposition method. Let 
$$V^m:=[x_{L}^{m},x_{R}^{m}],\quad x_L^1=x_L,\quad x_R^M=x_R,\quad m=1,2,\ldots,M$$
be a uniform partition of $D=[x_L,x_R]$ with the spatial step size $\Delta x$, thus $D=\cup_{m=1}^MV^m$. Further, we consider a uniform partition of $V^m$, $m=1,2,\ldots,M$,  with $J+1$ grid points in each fine interval, i.e.,
$$x_L^m=x_0^m<x_1^m<\cdots<x_{N-1}^m<x_{N}^m=x_R^m,\quad m=1,2,\ldots,M.$$
Differing from the traditional spatial partition, here we require the last two points of the element $V^m$ coincide with the first two points of the element $V^{m+1}$, that is 
$$x_{0}^{1}=x_L,\quad x_{J}^{M}=x_R,\quad x_{J-1}^{m}=x_{0}^{m+1},\quad x_{J}^{m}=x_{1}^{m+1},\quad m=1,2,\ldots,M-1.$$ 
In this situation, we remark that $\Delta x\neq(x_R-x_L)/J$. The general idea of the overlapping domain decomposition method is displayed in Fig. \ref{oddM_fig}.

After the above preliminaries, in element $V^m$, $m=1,2,\ldots,M$ we have
\begin{equation}\label{14}
\begin{split}
p^{m,n+1}=p^{m,n}+\frac{\tau}{2}\Big(q_{xx}^{m,n+1}+q_{xx}^{m,n}\Big),\quad q^{m,n+1}=q^{m,n}-\frac{\tau}{2}\Big(p_{xx}^{m,n+1}+p_{xx}^{m,n}\Big),
\end{split}
\end{equation}
where $u^m=p^{m}+{\rm i}q^{m}$ is the solution of \eqref{NLS_two1} over the $m$-th element. 

To discrete these equations with the boundary conditions \eqref{bc}, we mainly use the Chebyshev--Gauss--Lobatto quadrature points (see, e.g., \cite[Chapter 4]{RefBoyd}) in the interval $[-1,1]$ of the form $$\eta_{j}=\cos\Big(\frac{J-j}{J}\pi\Big),\quad j=0,1,\ldots,J.$$ 
Since \eqref{14} holds for the interval $[x_{L}^{m},x_{R}^{m}]$, in order to use the Chebyshev interpolation technique we need to introduce the following transformation
\begin{equation}\label{15}
\begin{split}
\eta^{m}:[x_{L}^{m},x_{R}^{m}]&\rightarrow[-1,1],\quad y^m\mapsto\frac{2}{x_{R}^{m}-x_{L}^{m}}y^m-\frac{x_{R}^{m}+x_{L}^{m}}{x_{R}^{m}-x_{L}^{m}},
\end{split}
\end{equation}
then after a straightforward calculation we arrive at 
$$
\Delta x=\frac{x_R-x_L}{M+\frac{(1-M)(1-\cos(\frac{\pi}{J}))}{2}}.
$$
\begin{example}
	If $J=2$, then 
	$$\Delta x=\frac{2(x_R-x_L)}{M+1}.$$
\end{example}
\begin{example}
	If $M=1$, then $\Delta x=x_R-x_L$. After the uniform partition of $D=[x_L,x_R]$ with $J+1$ grid points, we derive the overlapping domain decomposition method (or classical finite difference method) with the uniform spatial step size $(x_R-x_L)/J$.
\end{example}

Now, we present the interpolation $p^{m}(\eta^{m},t)$ and $q^{m}(\eta^{m},t)$ on the Chebyshev--Gauss--Lobatt collocation points $\eta^{m}\in[-1,1]$ given by
\begin{equation}\label{16}
\begin{split}
&p^{m}(\eta^{m},t)=\sum_{j=0}^{J}\widetilde{p}_{j}^{m}(t)\phi_{j}^{m}(\eta^{m}):=\big(\Phi^{m}(\eta^m)\big)^{\top}P^{m}(t),\\
&q^{m}(\eta^{m},t)=\sum_{j=0}^{J}\widetilde{q}_{j}^{m}(t)\phi_{j}^{m}(\eta^{m}):=\big(\Phi^{m}(\eta^m)
\big)^{\top}Q^{m}(t),
\end{split}
\end{equation}
where $\Phi^{m}=[\phi_{0}^{m},\cdots,\phi_{J}^{m}]^{\top}$ with $\phi_j$ being the Lagrangian interpolation function, $P^{m}=[\widetilde{p}_{0}^{m},\cdots,\widetilde{p}_{J}^{m}]^{\top}$, and $Q^{m}=[\widetilde{q}_{0}^{m},\cdots,\widetilde{q}_{J}^{m}]^{\top}$.

Taking the second-order partial derivatives of $p^m$ and $q^m$ with respect to $y^m$ and using the transformation \eqref{15}, it holds that
\begin{equation}\label{17}
\begin{split}
\frac{\partial^{(2)}p^{m}(y^{m},t)}{\partial y^{m,(2)}}=\Big(\frac{2}{x_{R}^{m}-x_{L}^{m}}\Big)^{2}\frac{\partial^{(2)}p^{m}(\eta^{m},t)}{\partial \eta^{m,(2)}}=\Big(\frac{2}{x_{R}^{m}-x_{L}^{m}}\Big)^{2}\big(\Phi^{m}\big)^\top D^{m,(2)}P^{m},\\[2mm]
\frac{\partial^{(2)}q^{m}(y^{m},t)}{\partial y^{m,(2)}}=\Big(\frac{2}{x_{R}^{m}-x_{L}^{m}}\Big)^{2}\frac{\partial^{(2)}q^{m}(\eta^{m},t)}{\partial \eta^{m,(2)}}=\Big(\frac{2}{x_{R}^{m}-x_{L}^{m}}\Big)^{2}\big(\Phi^{m}\big)^\top D^{m,(2)}Q^{m},
\end{split}
\end{equation}
respectively, where $D^{m,(2)}$ is a $(J+1)\times(J+1)$ differential matrix on the element $V^m$. It can be checked that the entries of $D^{m,(r)}$ for any $r\geq1$ are given by
\begin{equation}\label{18}
[D^{m,(r)}]_{i+1,j+1}=\left\{\begin{array}{cc}
\frac{r}{\eta_{j}^{m}-\eta_{i}^{m}}\Big(\frac{e_{j}}{e_{i}}[D^{m,(r-1)}]_{i+1,i+1}-[D^{m,(r-1)}]_{i+1,j+1}\Big),&~~i\neq j,\\[3mm]
-\sum\limits_{k=0,k\neq i}^{J}[D^{m,(r)}]_{i+1,k+1},&~~i=j
\end{array}\right.
\end{equation}
with $e_{0}=(-1)^{J}/2,~e_{J}=1/2$, $e_{j}=(-1)^{J-j},~j=1,2,\ldots,J-1$. We refer to \cite{RefBoyd,RefWelfert} for more details. Particularly, the entries of $D^{m,(1)}$ are defined as
\begin{equation*}
D^{m,(1)}_{i,j}=\left\{\begin{array}{cc}
-\frac{2N^{2}+1}{6},& \text{if}~~i=j=0,\\[2mm]
\frac{2N^{2}+1}{6},& \text{if}~~i=j=J,\\[2mm]
\frac{\eta_{i}^{m}}{2[1-(\eta_{i}^{m})^{2}]},& ~~~~\text{if}~~i=j\neq 0,~J,\\[4mm]
-\frac{c_{i}}{c_{j}}\frac{1}{\eta_{i}^{m}-\eta_{j}^{m}},& \text{if}~~i\neq j~~~~~~~~
\end{array}\right.
\end{equation*}
with $c_{0}=2(-1)^{J},~c_{J}=2$, $c_{j}=(-1)^{J-j},~j=1,2,\ldots,J-1$.

Based on the above overlapping idea, we introduce the partition of $D$ with the following grid points  $$x_L=x_{0}^{1}<x_{1}^1<\cdots<x_{J-1}^{1}<x_{1}^{2}<\cdots<x_{J-1}^{m-1}<x_{1}^{m}<\cdots<x_{J-1}^M<x_{J}^{M}=x_R.$$
Denote by
\begin{equation*}
\begin{split}
&\widetilde{\Phi}=[\phi_{0}^{1}(\eta^{1}),\cdots,\phi_{J-1}^{1}(\eta^{1}),\phi_{1}^{2}(\eta^{2}),\cdots,\phi_{J}^{M}(\eta^{M})]^{\top}\in\mathbb{R}^{M(J-1)+1},\\[2mm]
&\widetilde{P}=[\widetilde{p}_{0}^{1}(t),\cdots,\widetilde{p}_{J-1}^{1}(t),\widetilde{p}_{1}^{2}(t),\cdots,\widetilde{p}_{J}^{M}(t)]^{\top}\in\mathbb{R}^{M(J-1)+1},\\[2mm]
&\widetilde{Q}=[\widetilde{q}_{0}^{1}(t),\cdots,\widetilde{q}_{J-1}^{1}(t),\widetilde{q}_{1}^{2}(t),\cdots,\widetilde{q}_{J}^{M}(t)]^{\top}\in\mathbb{R}^{M(J-1)+1}.
\end{split}
\end{equation*}
It follows from \eqref{16} and the fact $D=\cup_{m=1}^MV^m$ that the interpolations of functions $p$ and $q$ in arbitrary collocation point $x\in D$ read as
\begin{equation*}
p(x)=\widetilde{\Phi}^\top\widetilde{P},\quad q(x)=\widetilde{\Phi}^\top\widetilde{Q},
\end{equation*}
which, along with \eqref{17} yields
\begin{equation}\label{20}
\frac{d^{(2)}p(x)}{dx^{(2)}}=\widetilde{\Phi}^\top\widetilde{D}^{(2)}\widetilde{P},\quad \frac{d^{(2)}q(x)}{dx^{(2)}}=\widetilde{\Phi}^\top\widetilde{D}^{(2)}\widetilde{Q},
\end{equation}
where
\begin{equation*}
\begin{split}
&\widetilde{D}^{(2)}_{[0:N-1,0:N]}=\frac{4}{\Delta x^{2}}D_{[0:J-1,0:J]}^{1,(2)},\\[2mm]
&\widetilde{D}^{(2)}_{[m(J-1)+1:(m+1)(J-1),m(J-1):(m+1)(J-1)+1]}=\frac{4}{\Delta x^{2}}D_{[1:J-1,0:J]}^{m+1,(2)},~~m=1,2,\ldots,M-2,\\[2mm]
&\widetilde{D}^{(2)}_{[(M-1)(J-1)+1:M(J-1)+1,(M-1)(J-1):M(J-1)+1]}=\frac{4}{\Delta x^{2}}D_{[1:J,0:J]}^{M,(2)},
\end{split}
\end{equation*}
and the remaining elements of the differential matrix $\widetilde{D}^{(2)}$ are zero.

Consequently, combining the temporal semi-discretization \eqref{first_order_system} and the Chebyshev interpolation \eqref{20}, we can obtain the following full discretization of \eqref{NLS_two1}:
\begin{align}\label{21}
U^{n+1}=U^n-\frac{{\rm i}}{2}\tau\Big(\widetilde{\Phi}^\top\widetilde{D}^{(2)}U^{n+1}+\widetilde{\Phi}^\top\widetilde{D}^{(2)}U^{n}\Big)
\end{align}
for all $n=0,1,\ldots,N-1$, where $U^{n}=\widetilde{P}^{n}+{\rm i}\widetilde{Q}^{n}$ and 
\begin{align*}
\widetilde{P}^n&=[\widetilde{p}_{0}^{1}(t_n),\cdots,\widetilde{p}_{J-1}^{1}(t_n),\widetilde{p}_{1}^{2}(t_n),\cdots,\widetilde{p}_{J}^{M}(t_n)]^{\top},\\[1.5mm] \widetilde{Q}^n&=[\widetilde{q}_{0}^{1}(t_n),\cdots,\widetilde{q}_{J-1}^{1}(t_n),\widetilde{q}_{1}^{2}(t_n),\cdots,\widetilde{q}_{J}^{M}(t_n)]^{\top}.
\end{align*}

\vspace{3mm}
\noindent{\bf (c). ODDS algorithm}
\vspace{2mm}

Based on the analytic expression \eqref{11} and the full discretization \eqref{21},  we have the following algorithm to compute the numerical solution to the one-dimensional stochastic NLS equation \eqref{NLS_one}.
\begin{algorithm}\label{a1}
	Choose the algorithm's parameters: time interval $[0,T]$; space domain $[x_L,x_R]$; temporal step size $\tau$; number of elements $M$; grid points $J$; orthonormal basis $\{e_k(x)\}_{k\geq1}$ and its truncation $\{e_k(x)\}_{k=1}^{K}$ to determine the $Q$-Wiener process $\Delta W_j^{K,n+1}$.
	
	Step 1. For each $n=1,2,\ldots,K-1$, $j=1,2,\ldots,M(J-1)$, take $u_{j}^{n}$ as the initial datum, solve \eqref{11} on the time interval $T_n$ and get
	\begin{equation*}\label{111}
	u_j^{\ast}=\exp\big\{-{\rm i}\big(\tau\lambda|u_j^{n}|^{2}+\varepsilon \Delta W_j^{n+1}\big)\big\}u_j^{n},
	\end{equation*}
	where $\Delta W_j^{n+1}=W(t_{n+1},x_j)-W(t_n,x_j)$.
	\vspace{2mm}
	
	Step 2. Let $U^{\ast}=(u_1^{\ast},u_2^{\ast},\cdots,u_{M(N-1)}^{\ast})^\top$. For each $n=1,2,\ldots,N-1$, take $U^{\ast}$ as the initial datum, solve \eqref{21} on the time interval $T_n$ and get
	\begin{equation*}\label{221}
	U^{n+1}=U^\ast-\frac{{\rm i}}{2}\tau\Big(\widetilde{\Phi}^\top\widetilde{D}^{(2)}U^{n+1}+\widetilde{\Phi}^\top\widetilde{D}^{(2)}U^{\ast}\Big).
	\end{equation*}
	
	Step 3. On the $n$-th time step (at time $t_n=n\tau$), generate the Gaussian random variables $\beta_k(t_n)$. According to \eqref{18}, compute the elements of $D^{m,(2)}$ for $m=1,2,\ldots,M$.
\end{algorithm}

\subsection{The ODDS algorithm for the multi-dimensional stochastic NLS equation}\label{2.2}
In this subsection, we present the ODDS algorithm for the $d$-dimensional stochastic NLS equation. Without loss of generality we restrict our discussion to the case $d=2$. Consider the following two-dimensional stochastic nonlinear system:
\begin{equation}\label{NLS_two}
{\rm i}du=\big[u_{xx}+u_{yy}+\lambda|u|^{2}u\big]dt+\varepsilon u\circ dW(t),\quad t\in(0,T]
\end{equation}
with an initial condition
$$u(0,x,y)=u_0(x,y),\quad x\in D=[x_L,x_R]\times[y_L,y_R]$$ 
and the boundary conditions 
\begin{equation*}\label{bcc}
\begin{split}
u(t,x_L,y)&=f_1(t),\quad u(t,x_R,y)=g_1(t),\quad{\rm on}\quad(0,T]\times\partial D,\\
u(t,x,y_L)&=f_2(t),\quad u(t,x,y_R)=g_2(t),\quad{\rm on}\quad(0,T]\times\partial D.
\end{split}
\end{equation*}

By using a similar technique as for the one-dimensional case, we split  \eqref{NLS_two} into the linear part and the nonlinear part. To define the ODDS algorithm for the two-dimensional case, the main difference lies in dealing with the linear part
\begin{align*}
{\rm i}du=(u_{xx}+u_{yy})dt.
\end{align*}
We use the local one dimensional idea to split the above equations as
\begin{align}
{\rm i}du&=u_{xx}dt,\label{1}\\[1.5mm]
{\rm i}du&=u_{yy}dt.\label{2}
\end{align}
Then, the algorithm developed in the previous subsection can be utilized to approximate the above four subsystems. The ODDS algorithm for the two-dimensional stochastic NLS equation \eqref{NLS_two} is presented as follows.
\begin{algorithm}\label{a2}
	Choose the algorithm's parameters: time interval $[0,T]$; space domain $[x_L,x_R]\times[y_L,y_R]$; temporal step size $\tau$; number of elements $M_1$ and $M_2$ in $x,y$-directions, respectively; grid points $J_1$ and $J_2$ in $x,y$-directions, respectively; orthonormal basis $\{e_{k_1,k_2}(x,y)\}_{k_1,k_2\geq1}$ and its truncation $\{e_k(x,y)\}_{k_1,k_2=1}^{K}$ to determine the $Q$-Wiener process $\Delta W_{j_1,j_2}^{K,n+1}$.
	
	Step 1. For each $n=1,2,\ldots,N-1$, $j_1=1,2,\ldots,M_1(J_1-1)$ and $j_2=1,2,\ldots,M_2(J_2-1)$, take $u_{j_1,j_2}^{n}$ as the initial datum, solve \eqref{11} on the time interval $T_n$ and get
	\begin{equation*}\label{1111}
	u_{j_1,j_2}^{\ast}=\exp\big\{-{\rm i}\big(\tau\lambda|u_{j_1,j_2}^{n}|^{2}+\varepsilon \Delta W_{j_1,j_2}^{n+1}\big)\big\}u_{j_1,j_2}^{n},
	\end{equation*}
	where $\Delta W_{j_1,j_2}^{n+1}=W(t_{n+1},x_{j_1},y_{j_2})-W(t_n,x_{j_1},y_{j_2})$.
	\vspace{2mm}
	
	Step 2. Let $U^{\ast}=(u_{1,1}^{\ast},u_{2,1}^{\ast},\cdots,u_{M_1(J_1-1),1}^{\ast},u_{1,2}^{\ast},u_{2,2}^{\ast},\cdots,u_{M_1(J_1-1),M_2(J_2-1)}^{\ast})^\top$. Take $U^{\ast}$ as the initial datum, solve \eqref{1} on the time interval $T_n$ by \eqref{21} and get
	\begin{equation*}\label{2211}
	U^{\ast\ast}=U^\ast-\frac{{\rm i}}{2}\tau\Big(\widetilde{\Phi}^\top\widetilde{D}^{(2)}U^{\ast\ast}+\widetilde{\Phi}^\top\widetilde{D}^{(2)}U^{\ast}\Big).
	\end{equation*}
	
	Step 3. For each $n=1,2,\ldots,N-1$, take $U^{\ast\ast}$ as the initial datum, solve \eqref{2} on the time interval $T_n$ by \eqref{21} and get
	\begin{equation*}\label{22111}
	U^{n+1}=U^\ast-\frac{{\rm i}}{2}\tau\Big(\widetilde{\Phi}^\top\widetilde{D}^{(2)}U^{n+1}+\widetilde{\Phi}^\top\widetilde{D}^{(2)}U^{\ast\ast}\Big).
	\end{equation*}
	
	Step 4. On the $n$-th time step (at time $t_n=n\tau$), generate the Gaussian random variables $\beta_{k_1,k_2}(t_n)$. According to \eqref{18}, compute the elements of $D^{m,(2)}$ for $m=1,2,\ldots,M_1$ and $m=1,2,\ldots,M_2$.
\end{algorithm}

\begin{remark}
	For the $S$-dimensional stochastic NLS equation with $S\geq3$, we only need to split the linear part of the considered system into $S$ subsystems
	$${\rm i}du=u_{x_sx_s}dt,\quad s=1,2,\ldots,S,$$
	and then use the similar algorithm as in Algorithm \ref{a2}.
\end{remark}
\section{Numerical experiments}\label{ne}
In this section we provide several numerical examples to illustrate the accuracy and capability of the algorithms developed in the previous section. We first present some preliminaries used throughout the following numerical implementation of Algorithms \ref{a1} and \ref{a2}.
\subsection{Preliminaries of the numerical implementation}
\begin{itemize}
	\item For $d=1$, we take the eigenvalues $\{\eta_k\}_{k=1}^{K}$ and the orthonormal basis $\{e_k\}_{k=1}^{K}$ of $L^{2}\big([x_L,x_R]\big)$ as
	$$e_k(x)=\sqrt{2}\sin(k\pi x),\quad \eta_k=1/k^3,$$
	which implies 
	\[
	\Delta W_{j}^{K,n+1}=\sum_{k=1}^{K}\sqrt{\frac{2}{x_{R}-x_{L}}}\sqrt{\frac{1}{k^3}}\sin\left(\frac{k\pi (x_j-x_{L})}{x_{R}-x_{L}}\right)\big(\beta_{k}(t_{n+1})-\beta_{k}(t_{n})\big).
	\]
	Here and in what follows, we take $K=500$.
	
	\item For $d=2$, we take the eigenvalues $\{\eta_{k_1,k_2}\}_{k_1,k_2=1}^{K}$ and the orthonormal basis $\{e_{k_1,k_2}\}_{k_1,k_2=1}^{K}$ of $L^{2}\big([x_L,x_R]\times[y_L,y_R]\big)$ as
	$$e_{k_1,k_2}(x,y)=2\sin(k_1\pi x)\sin(k_2\pi y),\quad \eta_{k_1,k_2}=1/(k_1^2+k_2^2)^2,$$
	which implies 
	\begin{align*}
	\Delta W_{j_1,j_2}^{K,n+1}=\sum_{k_1,k_2=1}^{K}&\frac{2}{k_1^2+k_2^2}\sqrt{\frac{1}{(x_{R}-x_{L})(y_{R}-y_{L})}}\sin\left(\frac{k_1\pi (x_{j_1}-x_{L})}{x_{R}-x_{L}}\right)\\[2mm]
	&\times\sin\left(\frac{k_2\pi (y_{j_2}-y_{L})}{y_{R}-y_{L}}\right)\big(\beta_{k_1,k_2}(t_{n+1})-\beta_{k_1,k_2}(t_{n})\big).
	\end{align*}
	\item Denote 
	$${\bf U}^n=[\widetilde{q}_{1}^{1}(t_n),\cdots,\widetilde{q}_{J-1}^1,\widetilde{q}_1^2,\cdots,\widetilde{q}_{J-1}^{M}(t_n),\widetilde{p}_{1}^{1}(t_n),\cdots,\widetilde{p}_{J-1}^1,\widetilde{p}_1^2,\cdots,\widetilde{p}_{J-1}^{M}(t_n)]^{\top},$$ then the full discretization \eqref{21} can be rewritten as an algebraic system:
	\begin{equation}\label{algebrain}
	(A\otimes B+C){\bf U}^{n+1}=(-A\otimes B+C){\bf U}^{n}+F,
	\end{equation}
	where $B=\widetilde{D}^{(2)}_{[1:M(J-1),1:M(J-1)]}$, $F\in\mathbb{R}^{2M(J-1)}$ describes the boundary conditions, and
	\begin{equation*}
	A=\left[\begin{array}{cc}
	-\frac{\tau}{2}&0\\[2mm]
	0&\frac{\tau}{2}
	\end{array}\right]_{2\times2},\quad C=\left[\begin{array}{cc}
	0&I\\[2mm]
	I&0
	\end{array}\right]_{2M(J-1)\times2M(J-1)}.
	\end{equation*}	
	We will compute \eqref{algebrain} using the Matlab command {\bf Algorithm 3.1}. Furthermore, once the differential matrix $D^{m,{2}}$ is known, then \eqref{algebrain} provides a feasible way to solve the stochastic NLS equation. The Matlab command {\bf Algorithm 3.2} relies on \eqref{18} to compute the elements of $D^{m,(2)}$ for $m=1,2,\ldots,M$. Since no confusion can arise, we simply drop the superscript $m,(2)$ on $D^{m,{2}}$.
	
	\begin{tabular*}{\columnwidth}{@{\extracolsep\fill}llllll@{\extracolsep\fill}}
		\toprule
		{\bf Algorithm 3.1} ~Code to compute the solution of a large sparse algebraic equation\\ $Gx=b$, where $G\in\mathbb{R}^{L\times L}$ is a sparse matrix and $b\in\mathbb{R}^{L}$ with $L=2M(J-1)$. The\\ input 
		$x0$ is an arbitrary non-zero column vector of length $L$.\\
		\midrule		
\begin{lstlisting}	%??????	
		
function x=matrix_solve(x0,G,b,L)

r=b-G*x0; u=zeros(length(x0),1);
while max(abs(r))>0.00001
 v(:,1)=r/norm(r); 
for j=1:m
 d=G*v(:,j);  
for i=1:j
 H(i,j)=v(:,i)'*d;
end
 u(:)=0;
for i=1:j
 u=H(i,j)*v(:,i)+u;
end
 u=d-u; H(j+1,j)=norm(u);
if (H(j+1,j)<0.0001||j==L)
 e=zeros(j+1,1); e(1)=norm(r);
 y=pinv(H(1:j+1,1:j))*e; 
 x0=x0+v(:,1:j)*y;r=b-G*x0;
 break;
end
v(:,j+1)=u/H(j+1,j);
end
end

\end{lstlisting}	\\
		\bottomrule
	\end{tabular*}
	
	\begin{tabular*}{\columnwidth}{@{\extracolsep\fill}llllll@{\extracolsep\fill}}
		\toprule
		{\bf Algorithm 3.2} ~Code to compute the differential matrix $D\in\mathbb{R}^{(J+1)\times(J+1)}$. \\
		\midrule		
\begin{lstlisting}	%??????	
		
function D=chebyshve_solve(J)

D=zeros(J+1,J+1);
K=(0:J)'; x=cos(pi*K/J); 
c=ones(J+1,1); c(1)=2; c(J+1)=2;

for k=1:J+1
 for j=1:J+1
  if (j==1&k==1)||(j==J+1&k==J+1)
   D(k,j)=(2*J^2+1)/6;
  elseif j==k
   D(k,j)=-x(k)/2/(1-x(k)^2);
  else
   D(k,j)=c(k)/c(j)*(-1)^(k+j)/(x(k)-x(j));
  end
 end
end

D(k,j)=-D(k,j); 

\end{lstlisting}	\\
		\bottomrule
	\end{tabular*}
	\item In order to demonstrate the efficiency and superiority of the proposed algorithm, we compare the ODDS algorithm with the following ones.
	
	(1) The stochastic multi-symplectic method (SMM for short; see \cite[Eq. (2.24)]{JWH2013}):
	\begin{equation}\label{25}
	\begin{split}
	{\rm i}\big(\delta_t^+u_{j+\frac12}^n+\delta_t^+u_{j-\frac12}^n\big)=&2\delta_x^+\delta_x^-u_j^{n+\frac12}+\lambda\big|u_{j+\frac12}^{n+\frac12}\big|^2u_{j+\frac12}^{n+\frac12}+\lambda\big|u_{j-\frac12}^{n+\frac12}\big|^2u_{j-\frac12}^{n+\frac12}\\[1.5mm]
	&+\varepsilon u_{j+\frac12}^{n+\frac12}\dot{W}_{j+\frac12}^{n+\frac12}+\varepsilon u_{j-\frac12}^{n+\frac12}\dot{W}_{j-\frac12}^{n+\frac12},
	\end{split}
	\end{equation}
	where $\delta_t^+u^n=(u^{n+1}-u^n)/\tau$, $\delta_x^+u_j=(u_{j+1}-u_j)/h_x$ and $\delta_x^-u_j=(u_{j}-u_{j-1})/h_x$. 
	
	(2) The finite difference splitting Crank--Nicolson scheme (FDSCN for short; see \cite[Eq. (57)]{CHLZ2019}):
	\begin{equation}\label{26}
	\begin{split}
	u_{j}^{\ast}&=u_j^{n}+{\rm i}\tau\Big(\delta_x^+\delta_x^-u_{j}^{n+\frac12,\ast}+\frac{\lambda}{2}\big(|u_j^n|^2+|u_{j}^{\ast}|^2\big)u_j^{n+\frac12,\ast}\Big),\\[1.5mm]
	u^{n+1}_j&=\exp\big(-{\rm i}\varepsilon\Delta W_j^{n+1}\big)u_j^{\ast},
	\end{split}
	\end{equation}
	where $u^{n+\frac12,\ast}=(u^n+u^\ast)/2$.
\end{itemize}

\subsection{Numerical examples}
After these preparations, now we concentrate on the numerical performance  of the ODDS algorithms.
\begin{example}
	{\rm In this example we show the soliton propagation at different instants of the following equation
		\begin{equation}\label{NLS_one1}
		{\rm i}du=\big[u_{xx}+|u|^{2}u\big]dt+\varepsilon u\circ dW(t),\quad t\in(0,T]
		\end{equation}
		in Figs. \ref{soliton1} and \ref{soliton2} with $\varepsilon=0.01$, $0.05$ and the initial condition
		$$u_0(x)=\sqrt{\frac{6}{5}}{\rm sech}\big(\sqrt{2}x\big)e^{{\rm i}x}.$$
		\begin{figure}[th!]
			\begin{center}
				\includegraphics[width=0.4\textwidth]{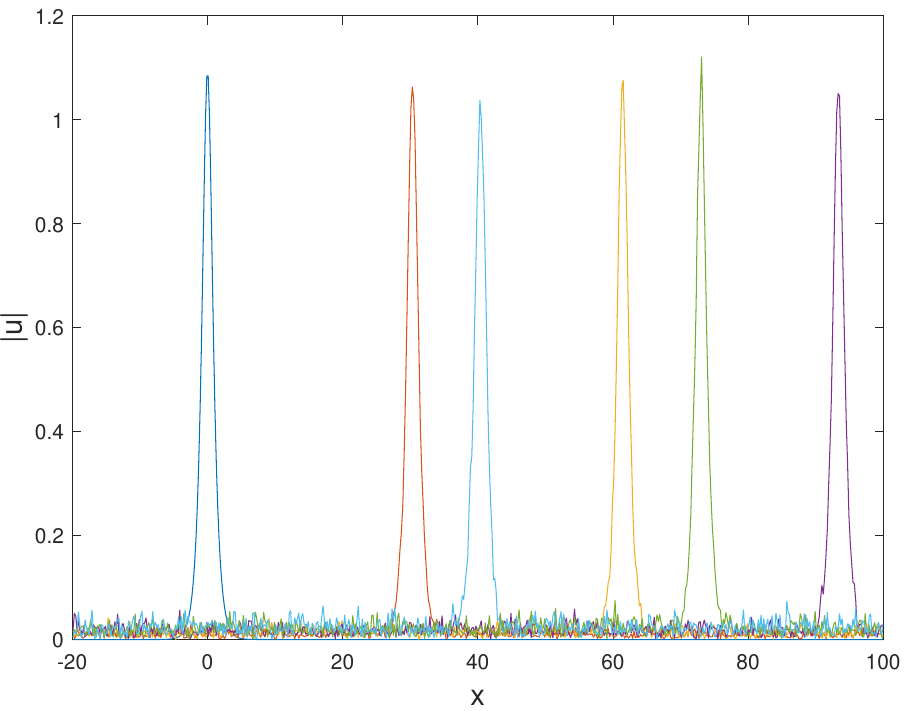}\qquad
				\includegraphics[width=0.45\textwidth]{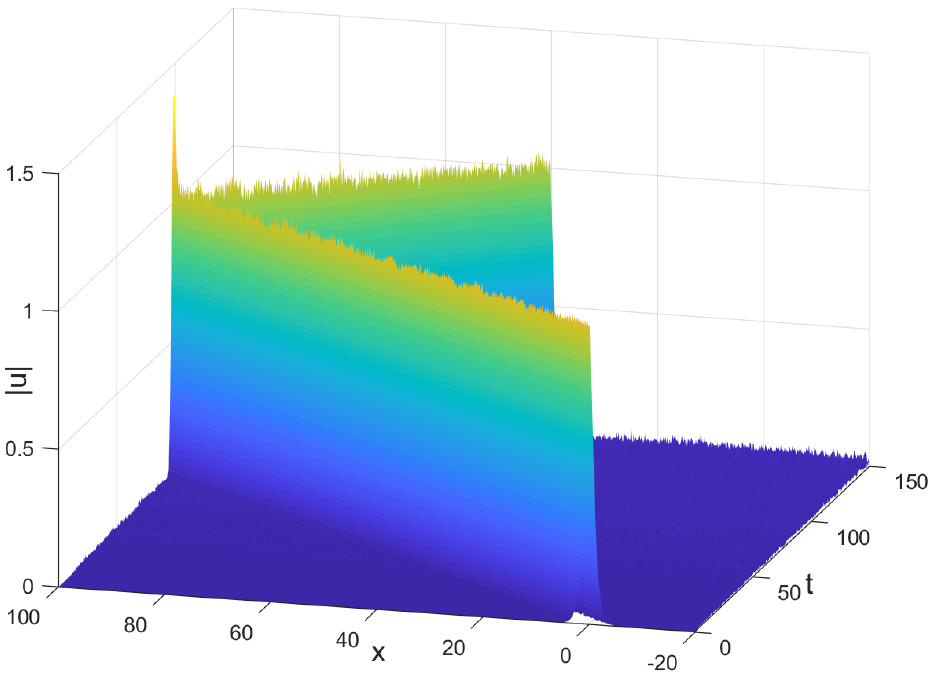}
				\caption{The soliton propagation of \eqref{NLS_one1} with the Dirichlet boundary conditions in $[-20,100]$ and $\varepsilon=0.01$. $J=30$, $M=10$, $T=150$, $\tau=0.015$. }\label{soliton1}
			\end{center}
		\end{figure}
		\begin{figure}[th!]
			\begin{center}
				\includegraphics[width=0.4\textwidth]{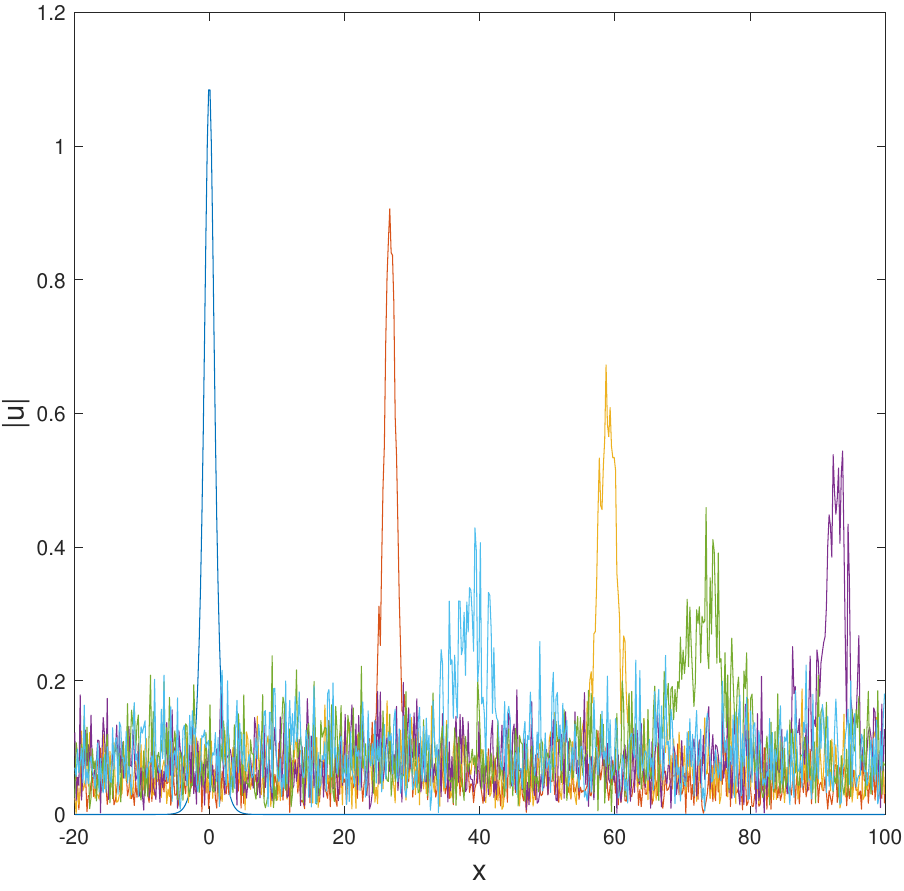}\qquad
				\includegraphics[width=0.45\textwidth]{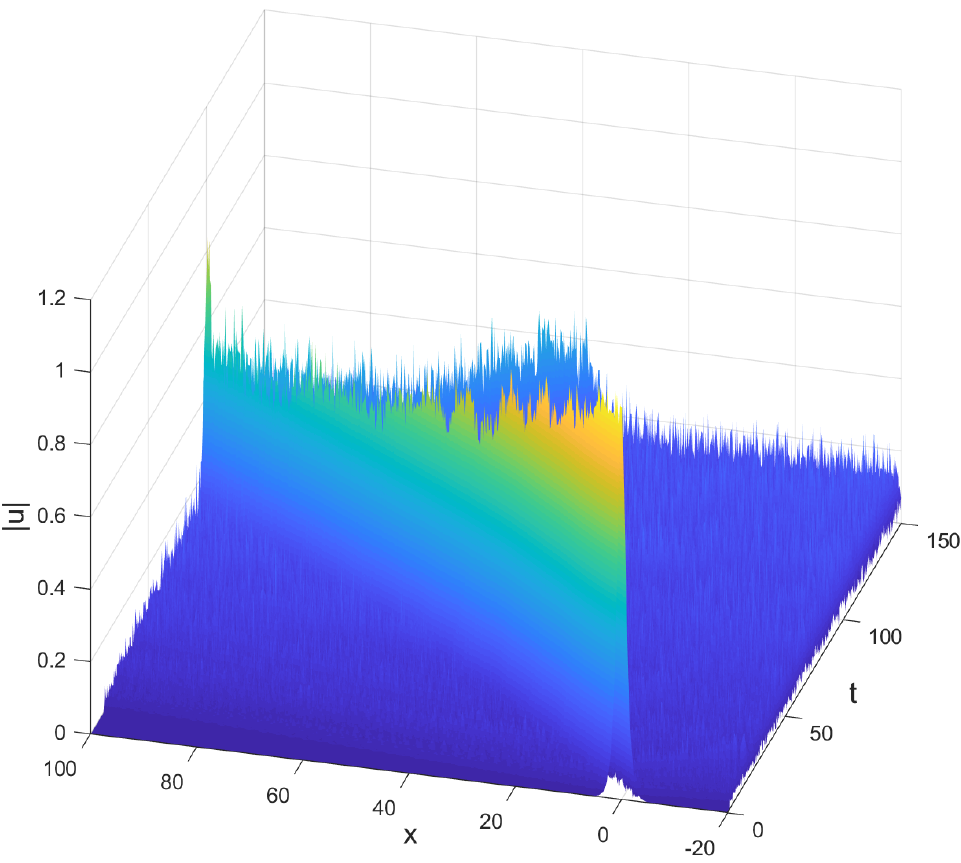}
				\caption{The soliton propagation of \eqref{NLS_one1} with the Dirichlet boundary conditions in $[-20,100]$ and $\varepsilon=0.05$. $J=30$, $M=10$, $T=150$, $\tau=0.015$.}\label{soliton2}
			\end{center}
		\end{figure}
		
		As is stated in \eqref{7}, equation \eqref{NLS_one1} possesses the charge conservation law almost surely under the homogeneous or periodic boundary conditions. Here, we verify this result by using our algorithm with zero Dirichlet boundary conditions. Fig. \ref{charge} presents the evolution of the discrete charge conservation law and the conservation errors of Algorithms \ref{a1} with $\varepsilon=0.01$ and $0.05$. We observe a good agreement with the continuous result.
		\begin{figure}[th!]
			\begin{center}
				\includegraphics[width=0.45\textwidth]{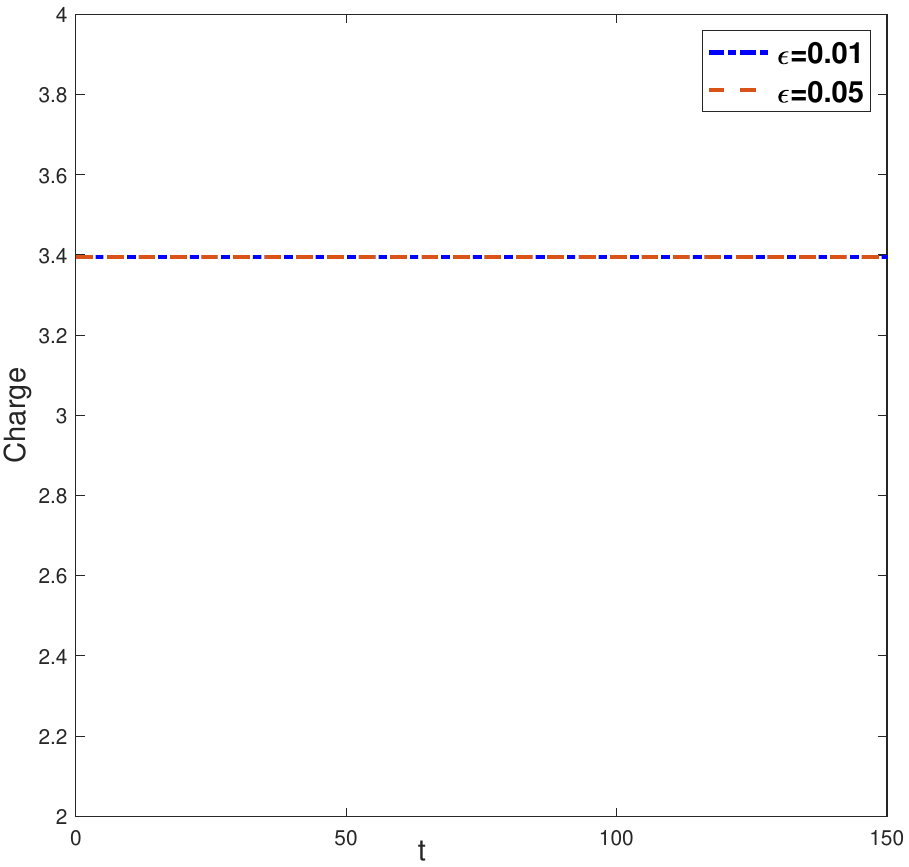}\qquad
				\includegraphics[width=0.45\textwidth]{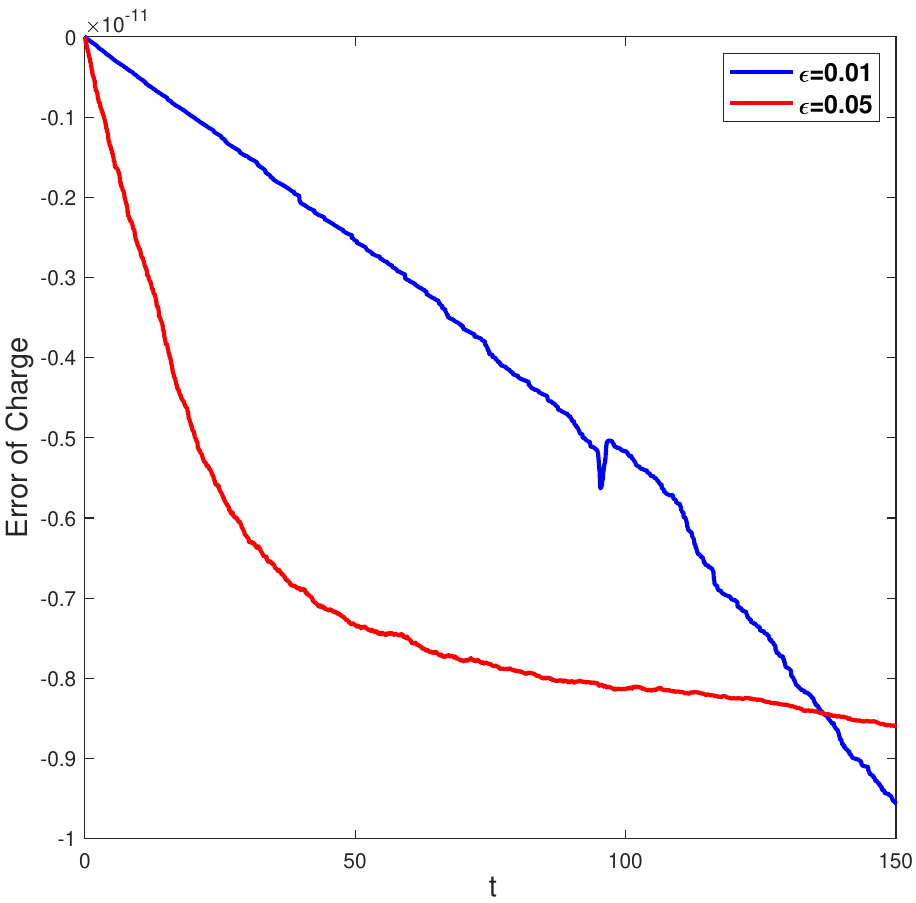}
				\caption{Evolution of the discrete charge (left), and the conservation error (right), over one trajectory with $\varepsilon=0.01, 0.05$. Zero Dirichlet boundary conditions in $[-20,100]$. $J=30$, $M=10$, $T=150$, $\tau=0.015$.}\label{charge}
			\end{center}
		\end{figure}
		
		Fig. \ref{energys} investigates the evolution of the discrete energy of the ODDS algorithm for different values of $\varepsilon=0.01$ and $0.05$, where the blues lines denote the discrete energies over 50 trajectories, the red line represents the discrete averaged energy, and the black line shows the discrete energy in the deterministic case. We see from the numerical experiment results that the discrete averaged energy possesses a linear growth property over 50 trajectories. 
		\begin{figure}[th!]
			\begin{center}
				\includegraphics[width=0.47\textwidth]{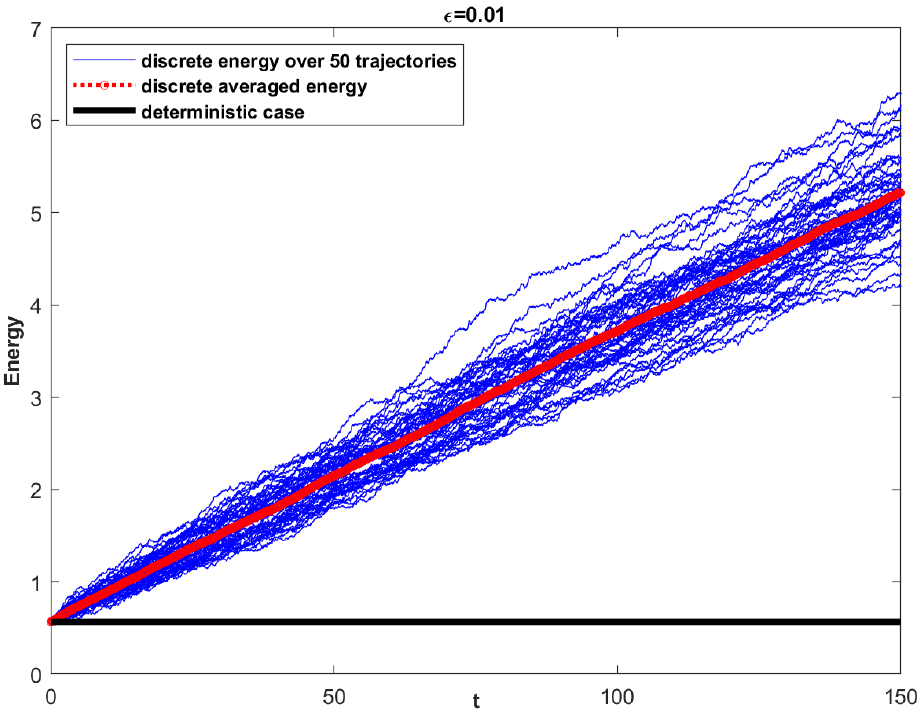}
				\includegraphics[width=0.47\textwidth]{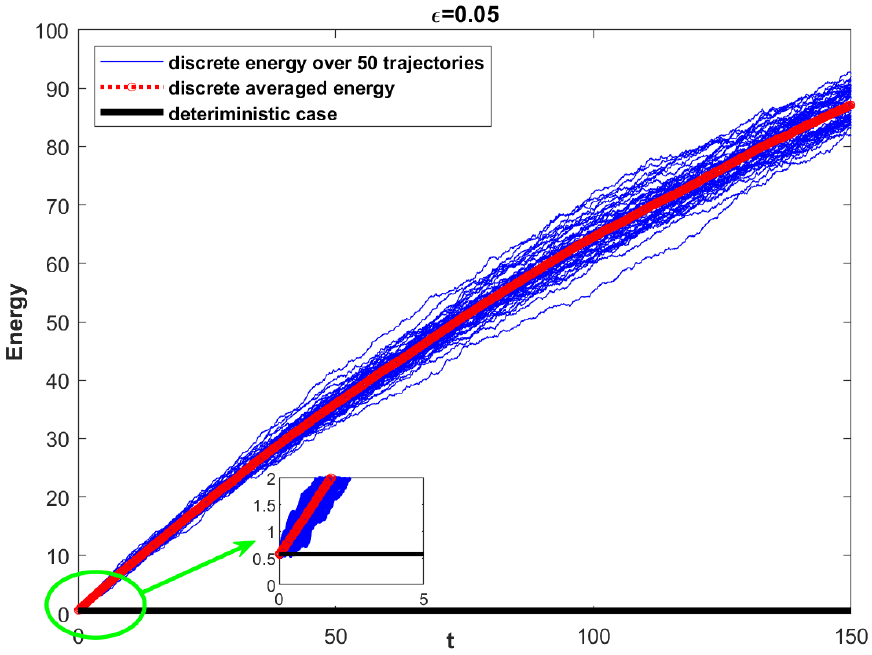}
				\caption{Evolution of the discrete energies for $\varepsilon=0.01$ (left), and $\varepsilon=0.05$ (right). Zero Dirichlet boundary conditions in $[-20,100]$. $J=30$, $M=10$, $T=150$, $\tau=0.015$.}\label{energys}
			\end{center}
		\end{figure}
		
		Now we compare the computational costs of the ODDS Algorithm \eqref{a1}, the SMM method \eqref{25} and the FDSCN scheme \eqref{26} for one-dimensional problem \eqref{NLS_one1} under the zero Dirichlet boundary conditions. Fig. \ref{cpu} demonstrates the computational efficiency of our ODDS algorithm in comparison with the SMM method and the FDSCN scheme. The reported CPU time is in seconds.
		\begin{figure}[th!]
			\begin{center}
				\includegraphics[width=0.65\textwidth]{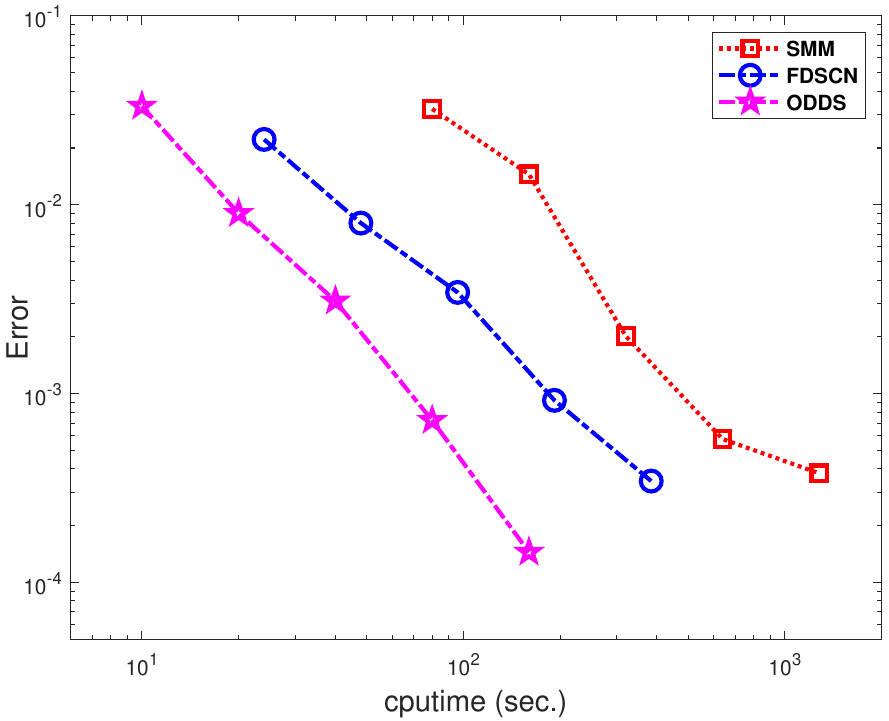}
				\caption{Efficiency for the ODDS algorithm, the SMM method and the FDSCN scheme for \eqref{NLS_one1}. Zero Dirichlet boundary conditions in $[-20,100]$. $J=30$, $M=20$, $T=150$, $\tau=0.015$. The mesh sizes of the SMM method and the FDSCN scheme are given by  $h_x=0.2$ (i.e., $J=600$, $M=1$).}\label{cpu}
			\end{center}
		\end{figure}
	}
\end{example}
\begin{example}
	{\rm In this example we present the double soliton collision of \eqref{NLS_one1} with the initial condition
		\begin{equation}\label{0}
		u_0(x)=\sqrt{\frac{6}{5}}{\rm sech}\big(\sqrt{2}x\big)e^{2ix}+\sqrt{\frac{3}{5}}{\rm sech}\big(\sqrt{2}(x-30)\big)e^{0.5i(x-30)}.
		\end{equation}
		The solution is simulated with the Dirichlet boundary conditions in $[-20,150]$ along one trajectory with $\varepsilon=0.01$.
		Figs. \ref{soliton_coll1}--\ref{soliton_coll3} show the double soliton collisions at different times $t=0$, $12$, $60$ for the real part $p$, the imaginary part $q$ and the module $u$, respectively. 
		\begin{figure}[th!]
			\begin{center}
				\includegraphics[width=0.45\textwidth]{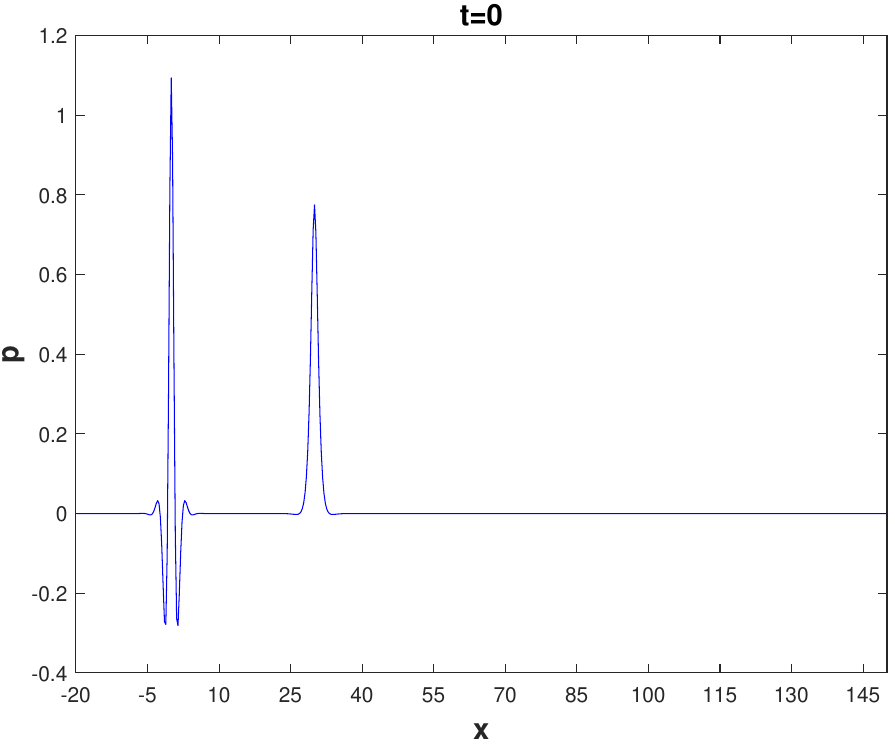}\qquad
				\includegraphics[width=0.45\textwidth]{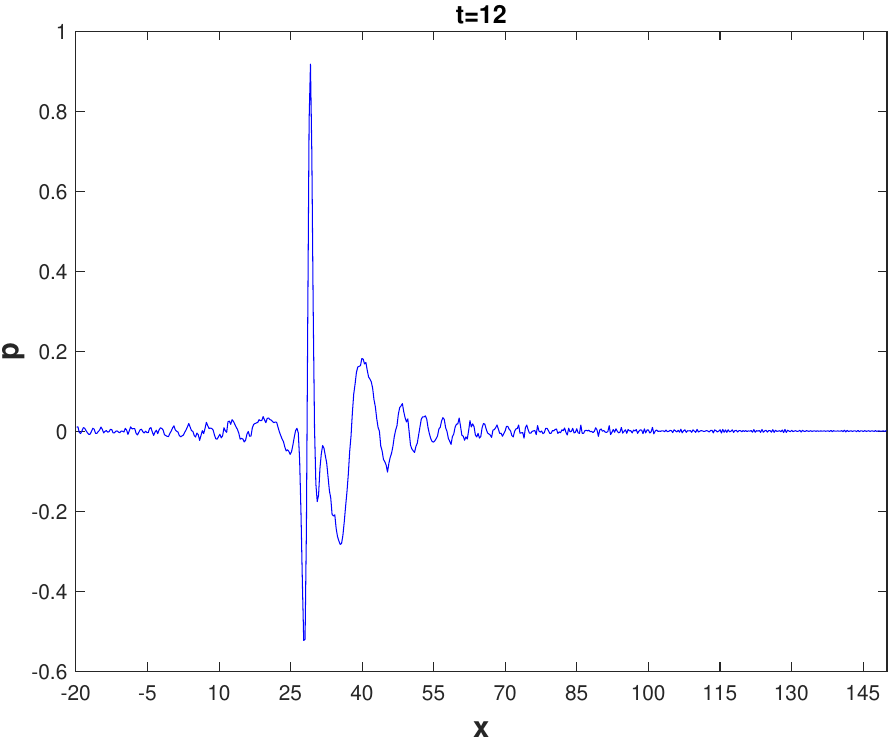}\\[1mm]
				\includegraphics[width=0.45\textwidth]{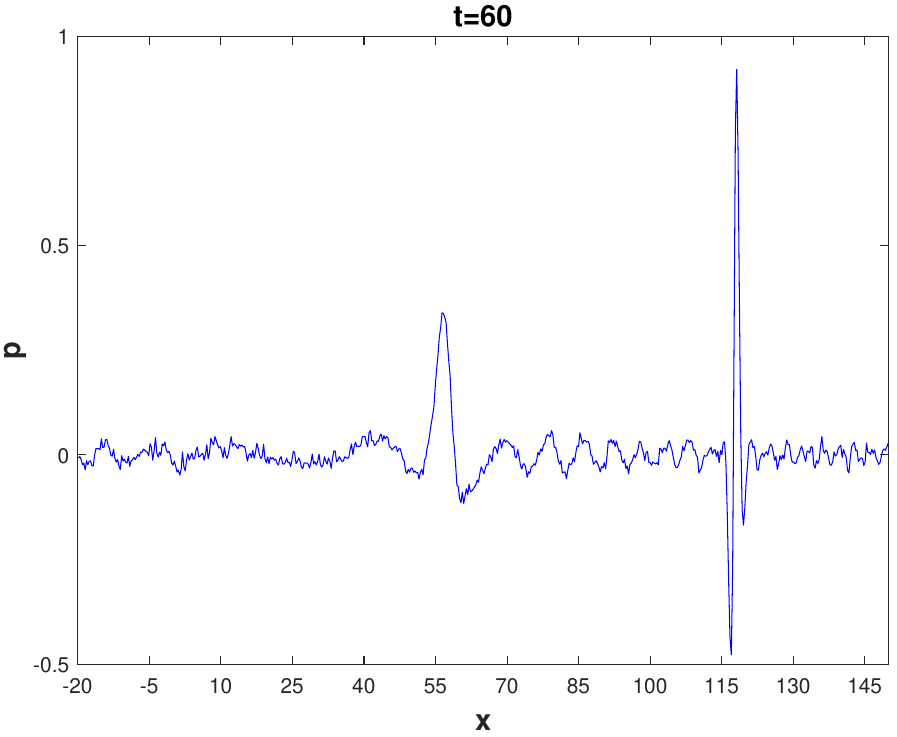}\qquad
				\includegraphics[width=0.46\textwidth]{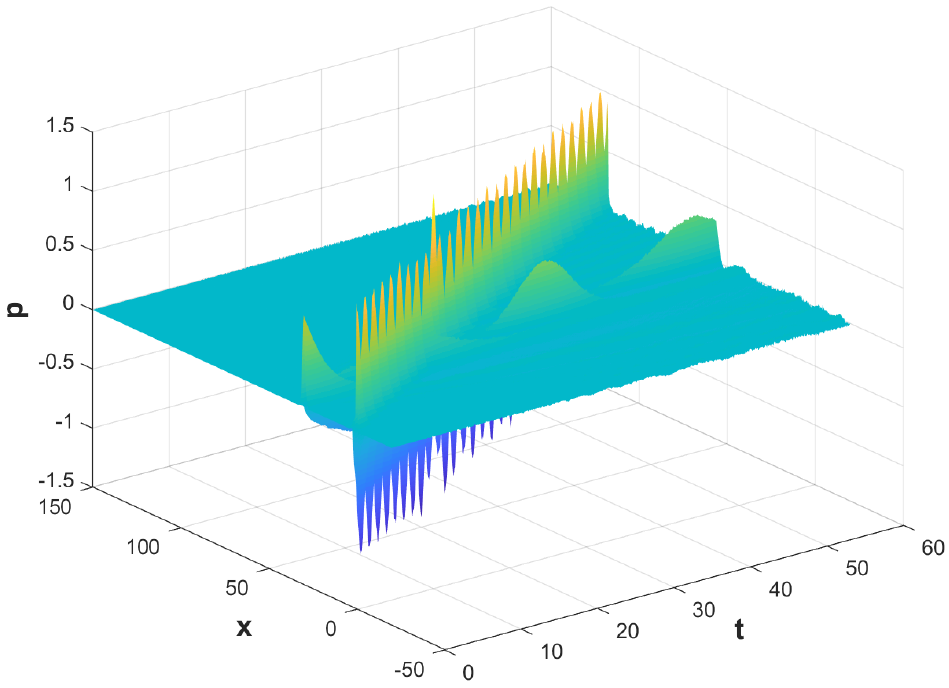}
				\caption{The double soliton collision of \eqref{NLS_one1} for the real part $p$ with the initial condition \eqref{0} along one trajectory. Dirichlet boundary conditions in $[-20,150]$. $J=20$, $M=5$, $T=60$, $\tau=0.006$.}\label{soliton_coll1}
			\end{center}
		\end{figure}	
		\begin{figure}[th!]
			\begin{center}
				\includegraphics[width=0.45\textwidth]{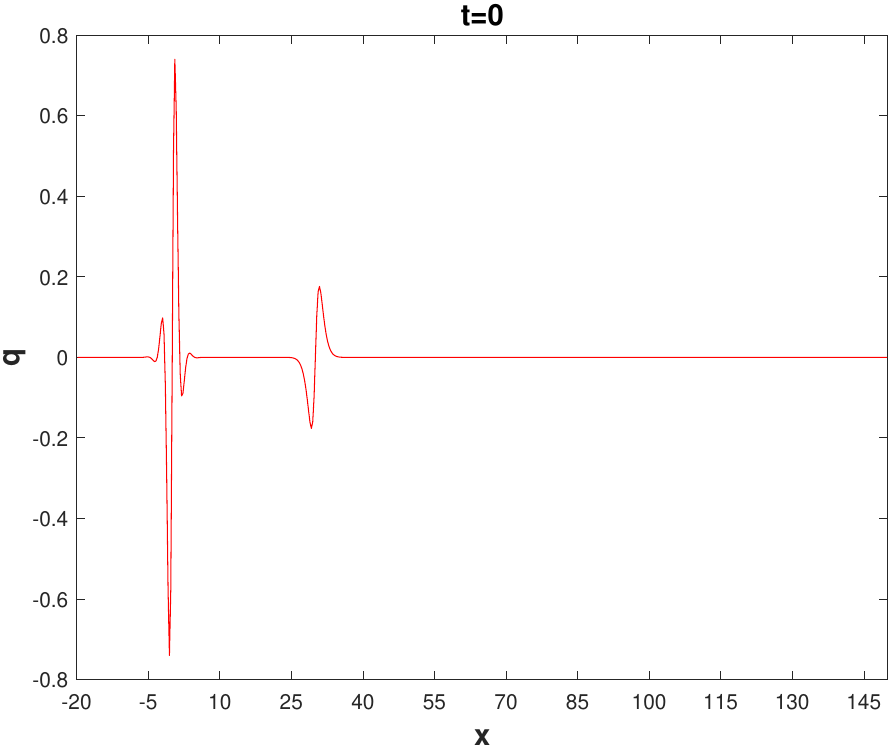}\qquad
				\includegraphics[width=0.45\textwidth]{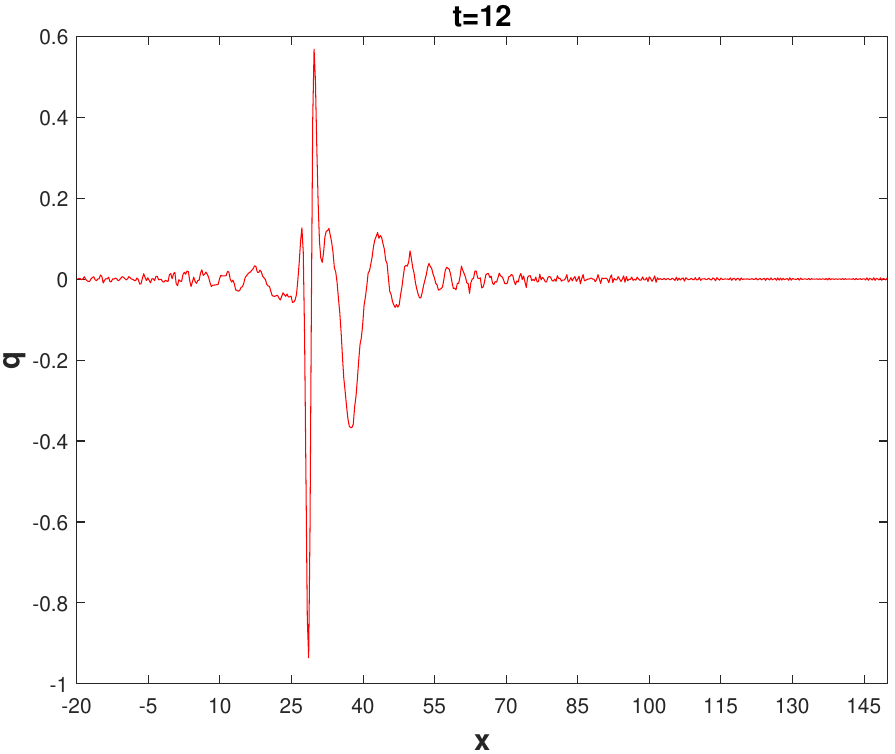}\\[1mm]
				\includegraphics[width=0.45\textwidth]{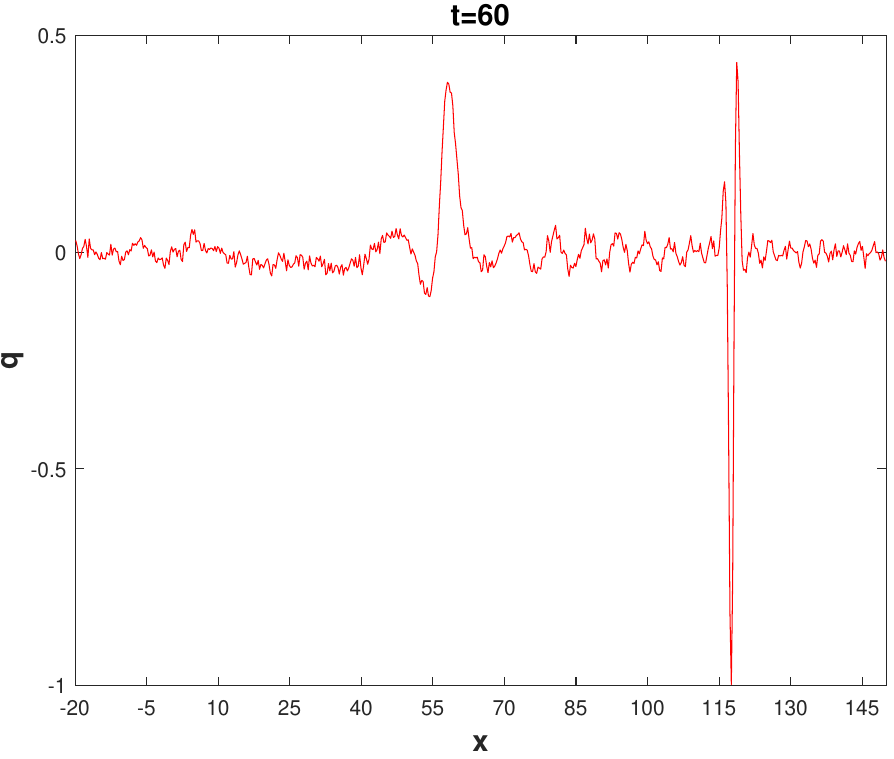}\qquad
				\includegraphics[width=0.46\textwidth]{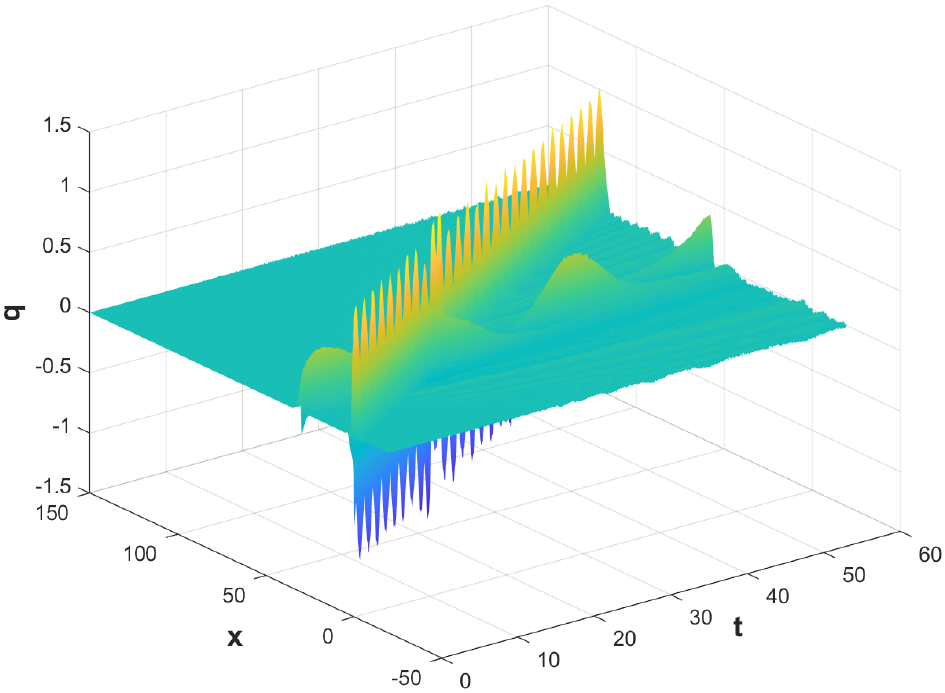}
				\caption{The double soliton collision of \eqref{NLS_one1} for the imaginary part $q$ with the initial condition \eqref{0} along one trajectory. Dirichlet boundary conditions in $[-20,150]$. $J=20$, $M=5$, $T=60$, $\tau=0.006$.}\label{soliton_coll2}
			\end{center}
		\end{figure}
		\begin{figure}[th!]
			\begin{center}
				\includegraphics[width=0.45\textwidth]{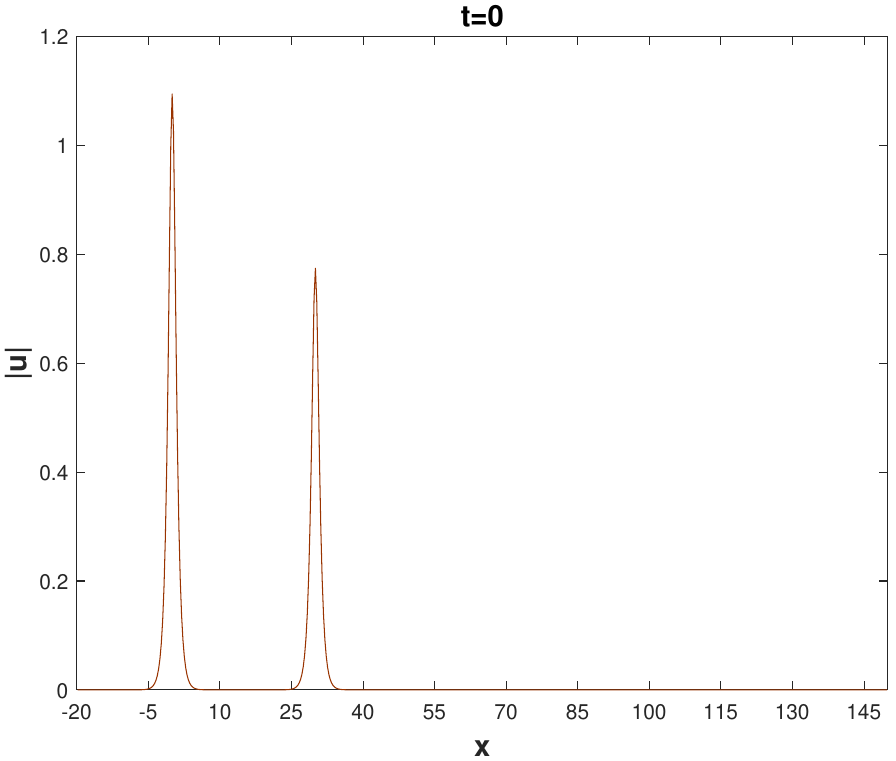}\qquad
				\includegraphics[width=0.45\textwidth]{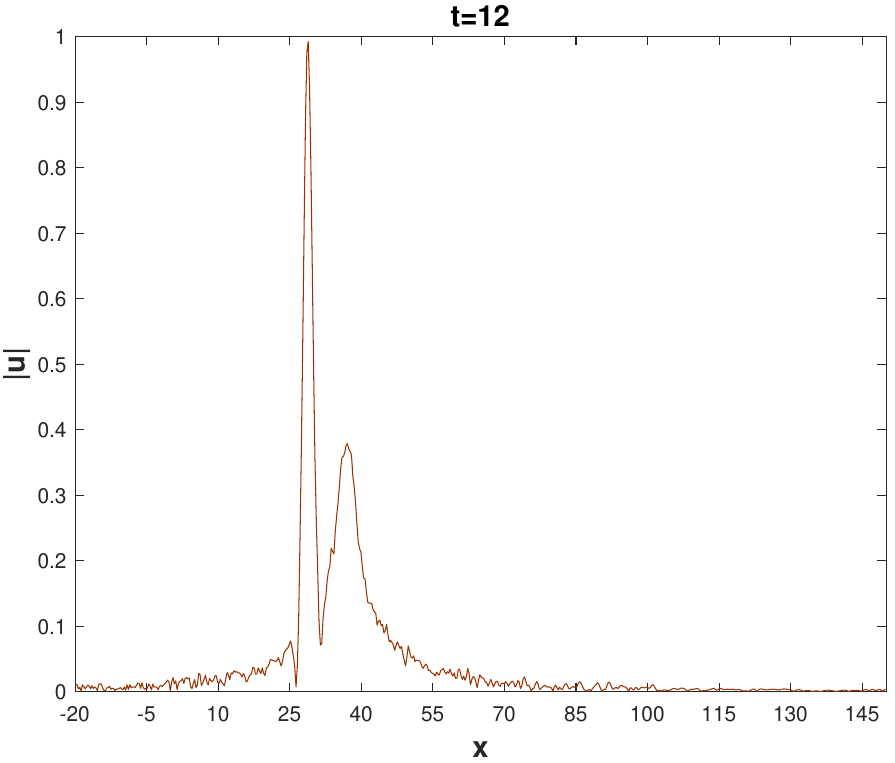}\\[1mm]
				\includegraphics[width=0.45\textwidth]{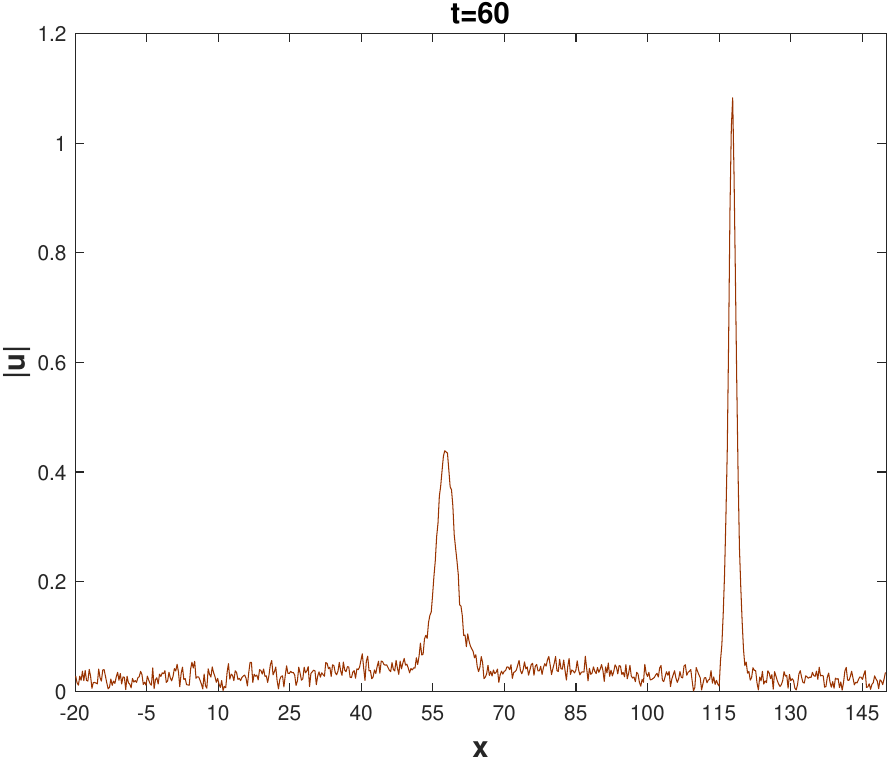}\qquad
				\includegraphics[width=0.46\textwidth]{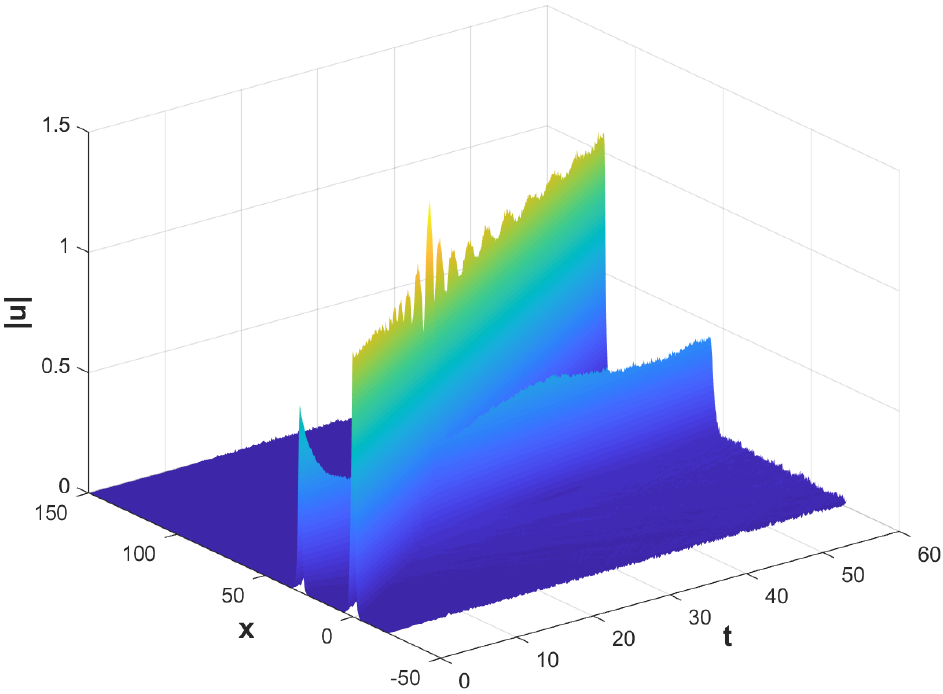}
				\caption{The double soliton collision of \eqref{NLS_one1} for the module $u$ with the initial condition \eqref{0} along one trajectory. Dirichlet boundary conditions in $[-20,150]$. $J=20$, $M=5$, $T=60$, $\tau=0.006$.}\label{soliton_coll3}
			\end{center}
		\end{figure}
	}
\end{example}
\begin{example}\label{ex3}
	{\rm In this example we consider the following two-dimensional stochastic NLS equation
		\begin{equation}\label{NLS_one2}
		{\rm i}du=\big[u_{xx}+u_{uu}+|u|^{2}u\big]dt+\varepsilon u\circ dW(t),\quad t\in(0,T].
		\end{equation}
		We choose the initial condition
		\begin{align}\label{00}
		u_0=A\exp{\big\{c_1x^2+c_2y^2\big\}},
		\end{align}
		where $A$, $c_1$ and $c_2$ are constants. The solution is computed with a Dirichlet boundary conditions in $[-10,10]\times[-10,10]$ with various sizes of the noise $\varepsilon=1$, $5$ and $10$. The results are presented in Figs. \ref{soliton_coll4}--\ref{soliton_coll6}.
		\begin{figure}[th!]
			\begin{center}
				\includegraphics[width=1\textwidth]{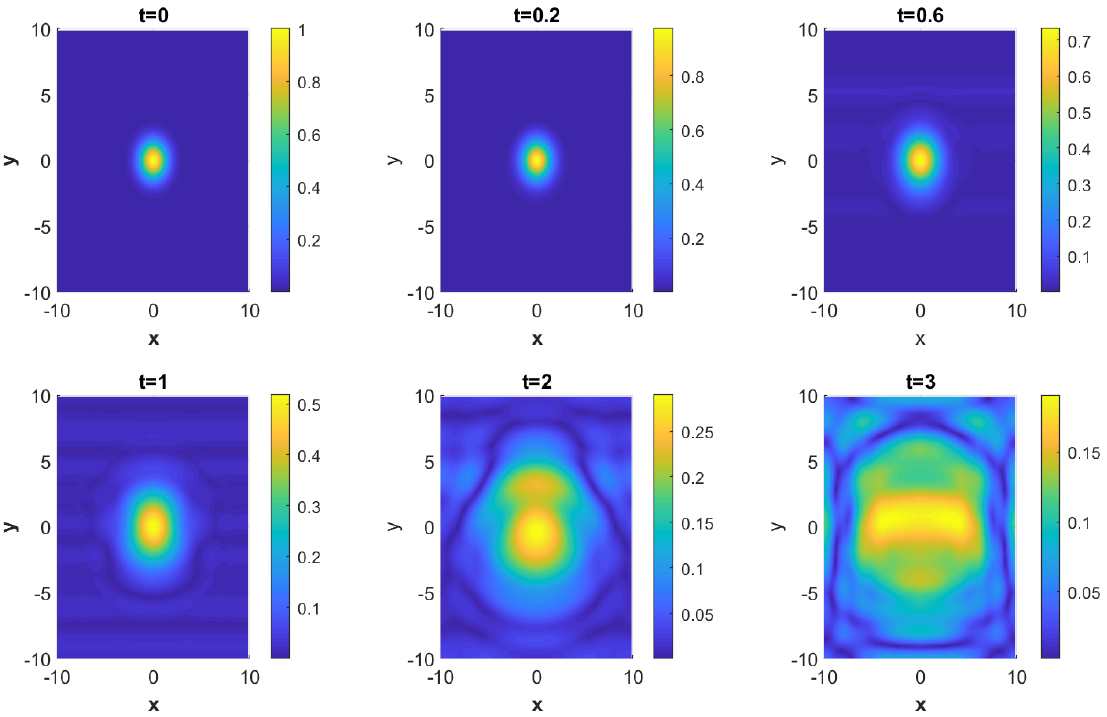}
				\caption{The solution of \eqref{NLS_one2} for the module $u$ with the initial condition \eqref{00} along one trajectory. Dirichlet boundary conditions in $[-10,10]\times[-10,10]$. $A=1$, $c_1=c_2=-1/2$, $J_1=J_2=32$, $M_1=M_2=4$, $T=3$, $\tau=0.01$, $\varepsilon=1$.}\label{soliton_coll4}
			\end{center}
		\end{figure}
		\begin{figure}[th!]
			\begin{center}
				\includegraphics[width=1\textwidth]{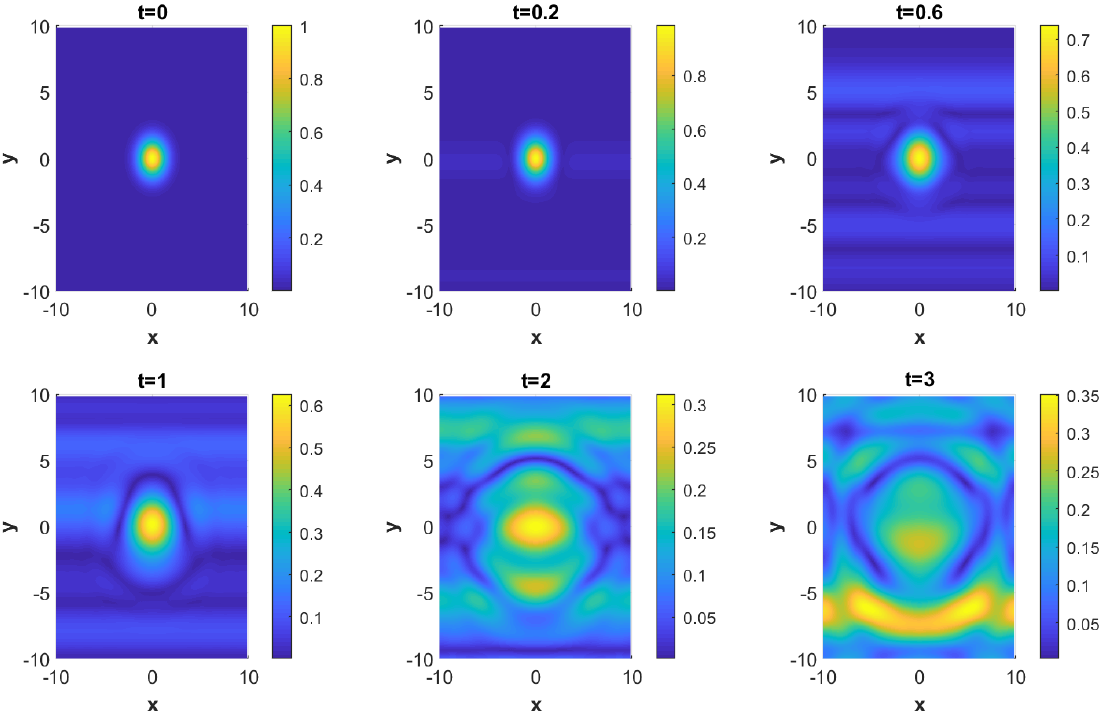}
				\caption{The solution of \eqref{NLS_one2} for the module $u$ with the initial condition \eqref{00} along one trajectory. Dirichlet boundary conditions in $[-10,10]\times[-10,10]$. $A=1$, $c_1=c_2=-1/2$, $J_1=J_2=32$, $M_1=M_2=4$, $T=3$, $\tau=0.01$, $\varepsilon=5$.}\label{soliton_coll5}
			\end{center}
		\end{figure}
		\begin{figure}[th!]
			\begin{center}
				\includegraphics[width=1\textwidth]{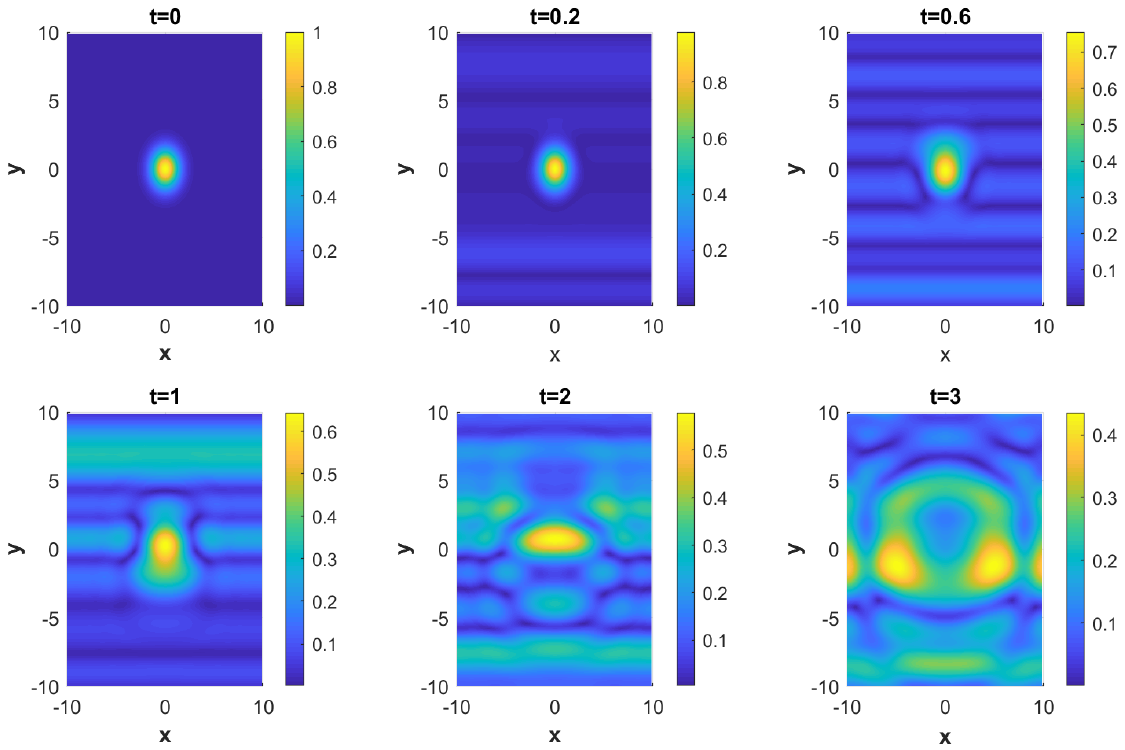}
				\caption{The solution of \eqref{NLS_one2} for the module $u$ with the initial condition \eqref{00} along one trajectory. Dirichlet boundary conditions in $[-10,10]\times[-10,10]$. $A=1$, $c_1=c_2=-1/2$, $J_1=J_2=32$, $M_1=M_2=4$, $T=3$, $\tau=0.01$, $\varepsilon=10$.}\label{soliton_coll6}
			\end{center}
		\end{figure}
	}
\end{example}

\begin{example}
	{\rm Without loss of generality, in this numerical example we restrict our discussion to the two-dimensional stochastic NLS equation \eqref{NLS_one2} with $\varepsilon=1$ to show the computational efficiency of the ODDS algorithm. We work in the same setting as in Example \ref{ex3}. 
		
		First, we apply the SMM method \eqref{25} and the FDSCN scheme \eqref{26} to the two-dimensional problem \eqref{NLS_one2}. Fig. \ref{cpu2} presents the computational cost of our ODDS algorithm in comparison with the SMM method and the FDSCN scheme. The reported CPU time is in seconds. From the figure, we can see that the ODDS algorithm can reduce the heavy computational load and is highly competitive.
		\begin{figure}[th!]
			\begin{center}
				\includegraphics[width=0.65\textwidth]{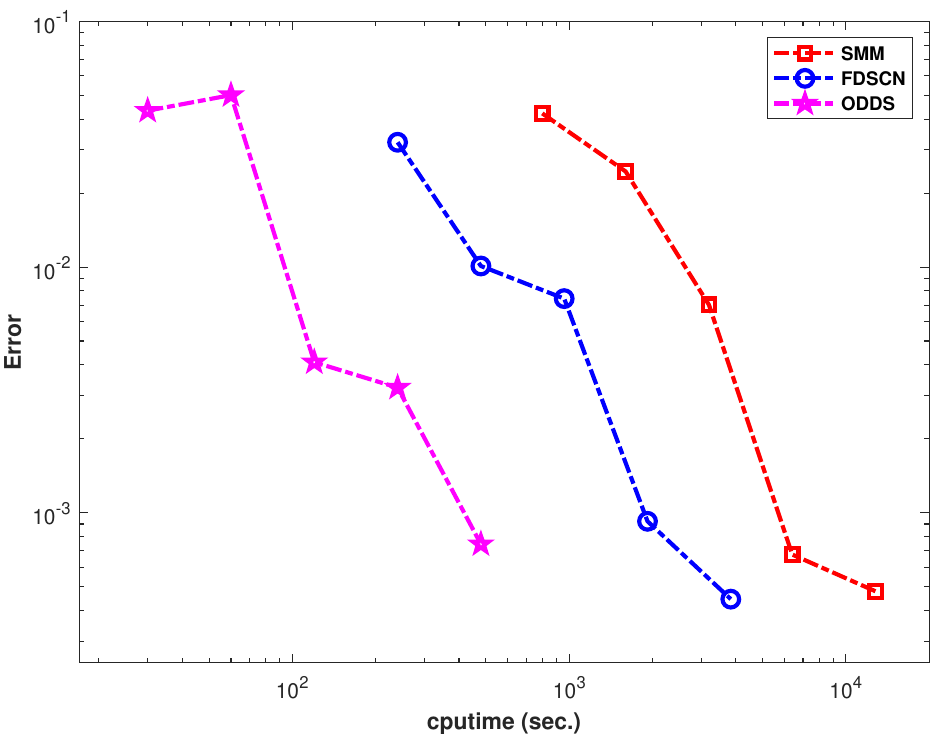}
				\caption{Efficiency for the ODDS algorithm, the SMM method and the FDSCN scheme for \eqref{NLS_one2} along one trajectory. Zero Dirichlet boundary conditions in $[-10,10]\times[-10,10]$. $A=1$, $c_1=c_2=-1/2$, $J_1=J_2=32$, $M_1=M_2=4$, $T=3$, $\tau=0.01$, $\varepsilon=1$. The mesh sizes of the SMM method and the FDSCN scheme are given by $h_x=h_y=5/32$ (i.e., $J_1=J_2=128$, $M_1=M_2=1$).}\label{cpu2}
			\end{center}
		\end{figure}
		
	}
\end{example}

\begin{example}
	{\rm We conclude this section with the mean-square convergence order in the temporal direction of the proposed ODDS algorithm for the one-dimensional stochastic NLS equation \eqref{NLS_one} with the initial condition $u_0=\sin(\pi x)$ 
		
		To compute the mean-square error, we run $P$ independent trajectories $u^{p}(t,\cdot)$ and $u^{p,n}(\cdot)$:
		$$
		Err:=\Big(\mathbb{E}\big[\|u(T,\cdot)-u^{N}(\cdot)\|_{l^2}^2\big]\Big)^{1/2}=\left(\frac{1}{P}\sum_{p=1}^P\big\|u^p(T,\cdot)-u^{p,N}(\cdot)\big\|_{l^2}^2\right)^{1/2}.
		$$
		We take time $T=1/4$, $[x_L,x_R]=[-1,1]$ and $P=500$. The reference solution is computed by the ODDS algorithm with small temporal step size $\tau=2^{-10}$. The number of trajectories $P=500$ is sufficiently large for the statistical errors not to significantly hinder the mean-square errors. The mean-square error is plotted in Table \ref{SchemeII_temporal}. The observed rates of convergence of the ODDS algorithm in time is close to 0.5$\sim$1. It is meaningful to give the mean-square convergence analysis theoretically in the future work.
		\begin{table}[th!]
			\caption{Mean-square errors of the ODDS algorithm for $\lambda=1$ and $\varepsilon=0.01$. \label{SchemeII_temporal}}%
			\begin{tabular*}{\columnwidth}{@{\extracolsep\fill}llllll@{\extracolsep\fill}}
				\toprule
				$\tau$&Err&Order\\
				\midrule
				$2^{-4}$  &5.1163E-1	 &--     \\
				$2^{-5}$  &2.6093E-1	&0.97\\
				$2^{-6}$  &1.2133E-1	&1.10\\
				$2^{-7}$  &7.0089E-2	&0.79\\	
				$2^{-8}$  &5.2949E-2	&0.40\\	
				$2^{-9}$  &2.7614E-2	&0.93\\	
				\bottomrule
			\end{tabular*}
		\end{table}
		
	}
\end{example}

\section{Concluding remarks}\label{cr}
The calculation of stochastic NLS equation is an interesting and important problem. One of the classical techniques is by the operator splitting. In this work, we have developed a high efficient ODDS algorithm to solve the stochastic NLS equation with a multiplicative noise by combining the splitting technique. Several numerical examples are presented to illustrate the capability of the algorithm. Although not considered in this work, this algorithm is flexible for the coupled stochastic NLS equation, the stochastic wave equation and the stochastic Maxwell equations, and has excellent computational efficiency.  One difficult and challenging future work is the mean-square convergence analysis of the ODDS algorithm.





\end{document}